\documentclass[aos,MSNbibl,seceqn,nameyear,dvips]{arximspdf}
\usepackage{graphicx}
\doi{10.1214/14-AOS1242} %kopijuoti is PTS
\volume{42}
\issue{4}
\pubyear{2014}
\firstpage{1657}
\lastpage{1688}
\docsubty{FLA}

\makeatletter
\def\arg{\operatorname{arg}}
\def\dim{\operatorname{dim}}
\newtheorem{Theorem}{Theorem}
\newproclaim{Example}{Example}
\newtheorem{Lemma}{Lemma}
\newtheorem{Proposition}[Lemma]{Proposition}
\newtheorem{Corollary}[Lemma]{Corollary}
\def\X{{\mathbf X}}
\def\H{{\mathbf H}}
\def\Z{\mathbf{Z}}
\def\c{\mathbf{c}}
\def\t{\mathbf{t}}
\def\u{\mathbf{u}}
\def\Cov{\operatorname{Cov}}
\def\sign{\operatorname{sign}}
\makeatother

\begin{document}
\begin{frontmatter}

\title{An adaptive composite quantile approach to~dimension reduction}
\runtitle{On adaptive quantile dimension reduction}

\begin{aug}
\author[a]{\fnms{Efang}~\snm{Kong}\corref{}\ead[label=e1]{e.kong@kent.ac.uk}}
\and
\author[b]{\fnms{Yingcun}~\snm{Xia}\thanksref{t1}\ead[label=e3]{staxyc@nus.edu.sg}}
\runauthor{E. Kong and Y. Xia}
\affiliation{University of Kent at Canterbury and National University
of Singapore}
\address[a]{School of Mathematics,\\
\quad Statistics and Actuarial Science\\
University of Kent at Canterbury\\
Kent\\
United Kingdom CT2 7NF\\
\printead{e1}} %adresu isvedimo komanda gale!
\address[b]{Department of Statistics\\\quad  and Applied Probability\\
National University of Singapore\\
Singapore 117546\\
\printead{e3}}
\end{aug}
\thankstext{t1}{Supported by National Natural Science Foundation of
China 71371095 and NUS Grant R-155-000-121-112.}

% HISTORY:
\received{\smonth{6} \syear{2013}}
\revised{\smonth{5} \syear{2014}}

% ABSTRACT
%
\begin{abstract}
Sufficient dimension reduction
[\textit{J.~Amer. Statist. Assoc.} \textbf{86} (1991) \mbox{316--342}]
has long been a
prominent issue in multivariate nonparametric regression analysis.
To uncover the central dimension reduction space, we
propose in this paper an adaptive composite quantile approach.
Compared to existing methods, (1) it requires minimal assumptions
and is capable of revealing all dimension reduction directions;
(2) it is robust against outliers and
(3) it is structure-adaptive, thus more efficient. Asymptotic results
are proved and numerical examples are provided, including a real
data analysis.
\end{abstract}

% KEYWORDS
\end{frontmatter}

%s1 #&#
\section{Introduction}\label{sec1}
Dimension reduction is a rather amorphous concept in statistics,
changing its characteristics and taking different forms
depending on the context.
In regression, the paradigm of sufficient dimension reduction
[\citet{r23}, \citeauthor{r6} (\citeyear{r6,r7})] which combines the idea of
dimension reduction with the concept of sufficiency, aims to
generate low-dimensional summary plot without appreciable loss of information.
In most cases, reductions are typically constrained to be linear
and the goal then is to estimate the central dimension reduction
subspace, or simply the central subspace.

\citet{r8} gave a formal definition and overviews
of the sufficient dimension reduction in regression,
which we adopt in this paper for the definition of the central subspace.
Suppose $ Y $ is a scalar dependent variable and
$ \X$ is the corresponding $p\times1$ vector of predictors. Let ${\mathbf B}$
be a ${p\times q}$ $(q\le p)$ (constant) orthonormal matrix and ${\mathbf
B}^\top$, its transpose. The space
${\mathcal S}({\mathbf B})$ spanned by the columns of ${\mathbf B}$,
is said to be
the (sufficient) dimension reduction subspace (DRS), if the
conditional distribution $F(\cdot|{\mathbf B}^\top\X)$ of $ Y $ given
${\mathbf B}^\top\X$ is identical to $ F(\cdot| \X)$, that is,
%
%e1.1 #&#
\begin{equation}
\label{drmodel} F(Y| \X) = F\bigl(Y|{\mathbf B}^\top\X\bigr) \qquad
\mbox{almost surely}.
\end{equation}
Consequently, a subspace is called a central subspace (CS), if it is not
only itself a DRS, but also a subset of any other DRS'. It thus
represents the minimal subspace that captured all the information
relevant to regressing $Y$ on $\X$. Under quite general
conditions, the CS exists and is given by
\[
{\mathcal S}_0 = \bigcap\bigl\{{\mathcal S}({\mathbf B})\dvtx
\mbox{model (\ref
{drmodel}) holds for } {\mathbf B}\bigr\};
\]
see \citet{r38} for the latest results on sufficient
conditions for the existence of CS.
Its dimension $\dim({\mathcal S}_0)=q(\le p)$ is referred to
as the structural dimension, while its orthogonal basis $\beta
_{01},\ldots,\beta_{0q}$ is called
the dimension reduction directions or simply the CS directions.
Let $ {\mathbf B}_0=(\beta_{01},\ldots,\beta_{0q})$, and thus
equivalent to
(\ref{drmodel}), we have
%
%e1.2 #&#
\begin{equation}
F(Y| \X) = F\bigl(Y|{\mathbf B}_0^\top\X\bigr) \qquad
\mbox{almost surely}. \label
{drmodel2}
\end{equation}
Research in dimension reduction methodologies, namely the search of
CS (directions), has garnered tremendous interest [\citet{r12}, \citet{r36}, \citet{r34},
\citet{r21}, \citet{r22}, \citet{r24}, \citet{r42} and \citet{r25}] since the seminal
work of \citet{r23}.
Some earlier research in this area such as \citet{r23},
was often based on either restrictive or hard-to-verify assumptions,
which limited their applications; while others
being model (moment)-based, targeted not at $S({\mathbf B}_0)$, but
instead
the reduction subspace $S({\mathbf B})$ associated with certain functional
of $F(Y| \X)$,
for example, the conditional mean [\citet{r9}] or the
conditional variance [\citet{r43}].
As we are going to demonstrate through the following example, such subspace
quite often is strictly a subset of~CS.
Consider the following model where
%
%e1.3 #&#
\begin{equation}
Y = \beta_1^\top\X+ \beta_2^\top\X
\varepsilon\quad\mbox{and}\quad E(\varepsilon|\X)=0.\label{ma3}
\end{equation}
As $E(Y|\X)= \beta_1^\top\X$, the central mean subspace $S(\beta
_1)$ is thus strictly contained in
$S(\beta_1,\beta_2)$, the full CS.

Seeing the restrictions with the aforementioned moment-based methods,
some consider the possibility of recovering all CS directions
by taking transformation of the response variable $Y $.
See, for example,
\citet{r42}, which practically requires assuming a parametric
model for
$\X$;
or \citet{r10}, where no theoretical results are available;
and \citet{r37}.
Others [\citet{r35}, \citet{r41}, \citet{r33}] tried to
extract information on CS
directly from the
conditional density or distribution function.
A~major drawback of the methodologies in the preceding four references
is that
the embedded estimation procedure is not structure-adaptive,
rendering the subsequent estimators of CS (directions) less
efficient. To see this, take model (\ref{ma3}), for example.
As the conditional density (distribution) function is nonlinear, the
smoothing parameter
used in constructing their
kernel estimators must therefore be small, that is, only a small
proportion of data is being used for local estimation. In contrary,
the conditional quantile function is in this case at least piecewise
linear, and consequently its estimation can be made more efficient
through the use of a larger (data-driven) bandwidth. Another reason
for us to consider a conditional quantile based approach is the
theoretical equivalence between conditional distribution functions
and conditional quantiles.

As in the case of conditional mean-based approach, we do not expect
the CS (directions) to be fully revealed via quantile regression at
any individual level. The solution we shall propose in this paper is
a combination of dimension reduction methods of \citet{r34} and
the composite quantile approach for regression [\citet{r44},
\citet{r13}, \citet{r11}], together with a
adaptive-weighting strategy.
The advantages of this new approach include:
(1) it requires minimal assumptions and can identify the CS
directions exhaustively; (2) it is robust against outliers, a
property inherited from quantile regression; and (3) the embedded
estimation procedure is structure-adaptive, that is, the use of a
data-driven bandwidth means more efficient use of data.

The paper is organized as follows. In Section~\ref{sec2}, we show how the CS
characterizes the composite outer product of gradients matrix. Based
on this characterization, Section~\ref{sec3} describes how an adaptive
composite quantile approach is \mbox{integrated} with the outer-product of
gradients (qOPG) method, and for comparison purposes, the composite
quantile minimum average variance method (qMAVE). In Section~\ref{sec4}, we
present regularity conditions and theoretical results on the
asymptotic normality of the qOPG estimator. Sections~\ref{sec5} and \ref{sec6} examine
some practical issues, such as bandwidth selection and determination
of the structural dimension. Section~\ref{sec7} contains some numerical
results, including an example of real data analysis. Section~\ref{sec8}
provides concluding remarks. All proofs are given in the \hyperref[app]{Appendix}.

%s2 #&#
\section{A composite quantile approach}\label{sec2}

Under model (\ref{drmodel2}), for any $0<\tau<1$, the $\tau$th
conditional quantile of $ Y $ given $ \X$,
\[
Q_\tau(\X)= \min\bigl\{y\dvtx F(y|\X) \ge\tau\bigr\}
\]
admits the following alternative expression:
%
%e2.1 #&#
\begin{equation}
Q_\tau(\X)= \min\bigl\{y\dvtx F\bigl(y|{\mathbf B}_0^\top
\X\bigr) \ge\tau\bigr\} = \tilde Q_\tau\bigl({\mathbf
B}_0^\top\X\bigr). \label{yuanquan}
\end{equation}
Its gradient vector
\[
\nabla Q_\tau(\mathbf{x})= \biggl[ \frac{\partial
Q_\tau(\mathbf{x})}{\partial x_1},\ldots,
\frac{\partial Q_\tau
(\mathbf{x})}{\partial
x_p} \biggr]^\top
\]
defined for any $\mathbf{x}=(x_1,\ldots,x_p)^\top\in R^p$, is thus related
to $\nabla\tilde Q_\tau(\cdot)$, the gradient vector of $\tilde
Q_\tau(\cdot)$, via the following identity:
%
%e2.2 #&#
\begin{equation}
\nabla Q_\tau(\mathbf{x}) = {\mathbf B}_0\nabla\tilde
Q_\tau\bigl({\mathbf B}_0^\top\mathbf{x}
\bigr).\label{kala}
\end{equation}
Consequently, we have the following fact for the corresponding
outer-product of gradients (OPG) matrix specific to level $\tau$:
%
%e2.3 #&#
\begin{eqnarray}\label{a2}
\Sigma(\tau)& =& E\bigl\{ \nabla Q_\tau(\X) \bigl[\nabla
Q_\tau(\X) \bigr]^\top\bigr\}
\nonumber\\[-8pt]\\[-8pt]
&=& {\mathbf B}_0E\bigl\{\nabla\tilde Q_\tau\bigl({
\mathbf B}_0^\top\X\bigr)\bigl[\nabla\tilde
Q_\tau\bigl({\mathbf B}_0^\top\X\bigr)
\bigr]^\top\bigr\}{\mathbf B}_0^\top.\nonumber
\end{eqnarray}
It is obvious that for any $\tau\in(0,1)$,
\[
{\mathcal S}\bigl(\Sigma(\tau) \bigr) \subseteq{\mathcal S}({\mathbf
B}_0).
\]
Indeed, plenty of examples exist where the above inequality holds
strictly for at least one $\tau\in(0,1)$. Consider, for example,
model (\ref{ma3}) with $\tau=0.5$ and the median of $\varepsilon$
equal to zero.
In other words, the CS may not be fully recovered by OPG matrices
specific to
any finite number of quantile levels.
The solution instead lies with the composite OPG matrix defined as
%
%e2.4 #&#
\begin{equation}
\Sigma= \int_0^1 \Sigma(\tau)\,d
\tau,\label{a5}
\end{equation}
as stated in the following lemma.

%
%le1 #&#
\begin{Lemma} \label{lemma1}
Suppose $\nabla Q_\tau(\cdot)$ exists for almost all $\tau\in(0,1)$
and $\X$.
We have
$
{\mathcal S}(\Sigma) = {\mathcal S}({\mathbf B}_0)$.
\end{Lemma}
%
%He et al (2013) observed a similar property of quantile for variable
%selection.

By definition, the composite OPG matrix $\Sigma$ is simply an
equally weighted average of the level-specific OPG matrices
$\Sigma(\tau)$, $0<\tau<1$. As previously demonstrated, $
\Sigma(\tau) $ for a given $\tau$ might contain little or no
information at all about the CS. Consider another example where $ Y
= \mathbf{x}_1 \varepsilon$, $ X = (\mathbf{x}_1, \ldots, \mathbf
{x}_p)^\top$ and
$\varepsilon$ has median zero. It is easy to see that
$\Sigma(0.5) = {\mathbf0}$, a $p\times p$ zero matrix.
We call such $ \Sigma(\tau) $ uninformative,
to which less weight should be assigned for the purpose of a more
revealing composite OPG matrix. Since whether or not any
level-specific
$ \Sigma(\cdot)$ is uninformative is not given a priori,
we suggest the following procedure to obtain an adaptively weighted
composite OPG matrix.
Suppose we have decided on the structural dimension $q$. For any given
$ \tau\in(0,1) $,
denote by $ \lambda_1(\tau)\ge\cdots \ge\lambda_p(\tau)\ge0 $, the
$p$ eigenvalues of $\Sigma(\tau)$.
The ``adaptively weighted'' composite OPG matrix is consequently
defined as
\[
\Sigma_w = \int_0^1 w(\tau)
\Sigma(\tau)\,d\tau,
\]
where the weight function
%
%e2.5 #&#
\begin{equation}
w(\tau) = \frac{\lambda_1(\tau)+\cdots + \lambda_q(\tau)}{\lambda
_1(\tau)+\cdots + \lambda_p(\tau)}, \label{weight}
\end{equation}
%
%{\color{red} Should we modify this as
%%\beginn
% w(\tau) = \frac{\lambda_1(\tau)+... +
reflects the percentage of information contained in the first $q $
eigenvectors of $\Sigma(\tau)$. If $ \Sigma(\tau) = {\mathbf0}$, we
define $ w(\tau) = 0 $.
Note that as ${\mathcal S}(\Sigma(\tau) ) \subseteq{\mathcal S}({\mathbf B}_0)$ for any $\tau$,
we have $ w(\tau) = 1 $ for any $\tau$ such that $\Sigma(\tau)>
{\mathbf0}$.
In practice, weights $w(\cdot)$ are derived from eigenvalues of estimates
of $\Sigma(\tau)$.

%s3 #&#
\section{Estimation of the dimension reduction directions}\label{sec3}
Based on Lemma \ref{lemma1}, the key to recovering the CS directions
lies with the estimation of the composite OPG matrix $ \Sigma$,
which in turn depends on the availability of a proper estimate of
the gradient vector $\nabla Q_\tau(\mathbf{x}) $ for any given $\tau
\in
(0,1)$ and $\mathbf{x}\in R^p$. Let $\hat\nabla Q_\tau(\mathbf{x})$\vadjust{\goodbreak}
denote such an
estimate. We can then construct estimate of the level-specific OPG
matrix~(\ref{a2}), and consequently estimate of the composite OPG
matrix~(\ref{a5}), as follows:
%
%e3.1 #&#
\begin{equation}
\hat\Sigma(\tau) = \frac{1}{n}\sum_{j=1}^n
\hat\nabla Q_\tau(\X_j),\qquad \hat\Sigma= \int
_0^1 \hat\Sigma(\tau)\,d\tau.\label{daoba}
\end{equation}
Various nonparametric estimators of $\nabla Q_\tau(\cdot) $ could be used
in (\ref{daoba}),
including kernel smoothing, nearest neighbor and spline
estimators; see, for example, \citet{r31}, \citet{r2}
and
\citet{r17}, \citet{r18}. In this paper, we opt for
the local polynomial estimators of \citet{r4} and \citet{r19}.
This is because, to show that $\hat\Sigma$ is root-$n$
consistent and asymptotically ``normal,'' we need the following two
prerequisites:
(i) $\hat\nabla Q_\tau(\mathbf{x})$ has a bias of order
$o_p(n^{-1/2})$ uniformly in $\mathbf{x}$ and in $\tau$;
(ii) a Bahadur-type expansion of $\hat\nabla Q_\tau(\mathbf{x})$,
again uniformly in $\mathbf{x}$ as well as in $\tau$.
Condition (i) can be met by approximating $Q_\tau(\cdot)$ locally with
polynomials in $p$ variables with high enough degrees. Condition
(ii), to be proved in the \hyperref[app]{Appendix} using results on empirical
processes and U-processes, extends what was obtained in \citet{r19},
where the uniformity is with respect to $\mathbf{x}$ only.

Suppose there exists some positive integer $k$ such that, for all
$\tau\in(0,1) $, $Q_\tau(\cdot)$ has partial derivatives of order up to
$k$ on $\mathcal{D}$, the compact support of $\X$ in $R^p$.
Consequently, for any given $\mathbf{x}=(x_1,\ldots,x_p)^\top\in
\mathcal{D}
$ and $\X$ near $\mathbf{x}$, $Q_\tau(\X)$ can be approximated by
its $k$th
order Taylor expansion, that is,
%
%e3.2 #&#
\begin{equation}
\label{mary} Q_\tau(\X)\approx Q_\tau(\mathbf{x})+\sum
_{1\le[u]\le k} \frac
{D^{\u} Q_\tau(\mathbf{x})}{\u!}(\X-\mathbf{x}%
)^{\u},
\end{equation}
where $\u=(u_{1},\ldots,u_{p})$ denotes a generic $p$-dimensional
vector of nonnegative integers, $[\u]=\sum_{i=1}^{p}u_{i}$,
$\u!=\prod_{i=1}^{p} u_{i}!$, $\mathbf{x}^{\u} = \prod_{i=1}^{p}
x_{i}^{ u_{i}}$ with the convention that $0^{0} = 1$, and $D^{\u}$
denotes the differential operator $\partial^{[\u]}/\partial x_{1}^{
u_{1}},\ldots, x_{p}^{ u_{p}}$. For ease of reference,
write $A=\{\u\dvtx  [\u]\le k\}$ and
$s(A) = \sharp(A)$, the cardinality of~$A$.

Suppose $(\X_i,Y_i)$, $i=1,\ldots,n$, are i.i.d. copies of $(\X,Y)$,
and $h_n$ is a smoothing parameter such that $h_{n}\to0$, as
$n\to\infty$. For any given $\mathbf{x}\in R^p$ and $\tau\in(0,1)$,
define two $s(A)\times1$ vectors as follows:
\begin{eqnarray*}
\mathbf{x}(h _{n},A)&=&\bigl(\mathbf{x}(h_{n},\u)
\bigr)_{\u\in A} \qquad\mbox{with } \mathbf{x}(h_{n},
\u)=h_{n}^{-[\u]}\mathbf{x} ^{\u},
\\
\c_{n}(\mathbf{x};\tau)&=&\bigl(c_{n,\u}(\mathbf{x};\tau )
\bigr)_{\u\in A}\qquad\mbox{with } c_{n,\u}(\mathbf{x};
\tau)=h_{n} ^{[\u]}D^{\u}Q_\tau(
\mathbf{x})/{\u!}.
\end{eqnarray*}
The local polynomial estimate
of $\c_{n}(\mathbf{x};\tau)$ is defined as a solution to the
following problem:
%
%e3.3 #&#
\begin{equation}
\min_{\c} \sum_{i=1}^{n}
\rho_{\tau}\bigl\{Y_{i}-\c^\top\X
_{i\mathbf{x}}(h_{n},A)\bigr\}K_{h_n}\bigl(|
\X_{i\mathbf{x}}|\bigr),\label{a8}
\end{equation}
where $\c= (c_{\u})_{\u\in A}\in R^{s(A)}$, $ \rho_\tau(s) = |s|
+ (2\tau-1 )s $, $\X_{i\mathbf{x}}=\X_i-\mathbf{x}$,
$|\cdot|$ stands for the supremum norm, $ K(\cdot) $ is a kernel
function in $R^p$ with finite support, and $ K_{h_n}(\cdot) =
K(\cdot/h_n)/h_n$.
Note that although in this paper we take $K(\cdot)$ to be the uniform
density function
on $[-1,1]^p$, the
$p$-dimensional cube in $R^p$, the results we obtain apply to
other cases such as the Epanechnikov kernel as well.

Since $\rho_\tau(s)\to\infty$, as $|s|\to\infty$,
solution to (\ref{a8}) always exists as long as $ K_{h_n}(|\X
_{i\mathbf{x}}|) > 0 $ for at least one $\X_i$.
Denote by $\hat\c_n(\mathbf{x};\tau)=(\hat c_{n,\u}(\mathbf
{x};\tau
))_{\u\in A}$, a solution to (\ref{a8}) and by $\hat\nabla Q_\tau
(\mathbf{x})$,
the local polynomial estimate of the gradient vector $\nabla Q_\tau
(\mathbf{x}) $:
\[
\hat\nabla Q_\tau(\mathbf{x}) =h_n^{-1}\bigl(
\hat\c_{n,\u}(\mathbf {x};\tau )\bigr)_{\u\in A,[\u]=1}.
\]
Consequently, we can construct estimates of the level-specific OPG matrix
$\Sigma(\tau)$ and of the composite OPG matrix $ \Sigma$ as follows:
%
%e3.4 #&#
\begin{eqnarray}\label{jon}
\hat\Sigma(\tau)&=&\frac{1}{n}\sum_{j=1}^n
\hat\nabla Q_\tau(\X_j)\bigl\{\hat\nabla Q_\tau(
\X_j)\bigr\}^\top;
\nonumber\\[-8pt]\\[-8pt]
\hat\Sigma &=& \int
_{0}^{1} \hat\Sigma(\tau) \,d\tau.\nonumber
\end{eqnarray}
For the sake of technical convenience, we focus on rather than the
$ \hat\Sigma$ in (\ref{jon}) but instead the following truncated
version:
%
%e3.5 #&#
\begin{equation}
\hat\Sigma_{\mathrm{T}} = \int_{\delta^*}^{1-\delta^*} \hat
\Sigma (\tau) \,d\tau, \label{a6}
\end{equation}
for some small $ \delta^* \in(0,1)$. This is due to the fact that
the uniformity in $\tau$ of the strong Bahadur-type representation of
$\hat\nabla Q_\tau(\mathbf{x})$ requires the conditional density of $Y$
given $\X$ at $Q_\tau(\X)$ to be uniformly bounded away from zero,
a condition apparently cannot be met by all $\tau\in(0,1)$. See
Lemma \ref{T1} and its proof given in the \hyperref[app]{Appendix} for more details.
Nevertheless, such truncation need not cause much concern. The
reasons are two-fold. On one hand, the
integral in (\ref{jon}) is approximated as
a summation over a sequence of
discretized $\tau$ values. On the other hand, the CS which is
derived from $\{Q_{\tau}(\cdot|\mathbf{x})\dvtx  0<\tau<1, \mathbf{x}\in
\mathcal{D}\} $ or
equivalently from $\Sigma$, is expected to closely resemble, if not
completely identical to, that from $\{Q_{\tau}(\cdot|\mathbf{x})\dvtx
\delta^*\le\tau\le1-\delta^*,\break  \mathbf{x}\in\mathcal{D}\} $ or
equivalently from
\[
\Sigma_{\mathrm{T}}=\int_{\delta^*}^{1-\delta^*} \Sigma(
\tau) \,d\tau,
\]
provided that $\delta^*>0$ is small enough. We assume this is indeed
the case, that is, $\Sigma_{\mathrm{T}}=\Sigma$.

As suggested at the end of Section~\ref{sec2}, we could further construct an
estimate of the adaptively-weighted truncated composite OPG matrix
as
%
%e3.6 #&#
\begin{equation}
\hat\Sigma_{w{\mathrm{T}}} = \int_{\delta^*}^{1-\delta^*} \hat
\Sigma (\tau) \hat w(\tau) \,d\tau, \label{a6A}
\end{equation}
with weight $ \hat w(\tau) $ calculated according to
formula (\ref{weight}) using the eigenvalues of $ \hat\Sigma(\tau)$.
However, to make sure less weights are assigned to those
uninformative matrices $\hat\Sigma(\tau)$ which are close to but not
exactly zero, we set $ \hat w(\tau) = 0 $ if the largest eigenvalue
of $\hat\Sigma(\tau)$ is below certain threshold.

In the ideal case where the structural dimension $q$ is known a
priori, estimates of the CS directions are simply given by the first
$q$ eigenvectors of $ \hat\Sigma_{\mathrm{T}}$: $\hat\beta_k$,
$k=1,\ldots,q$. Details on how to estimate $q$ when it is unknown as
well as bandwidth selection
are given in Sections~\ref{sec5} and \ref{sec6}, respectively. Similar to \citet{r34},
the above estimator can be further refined as follows.
Relabel the above obtained estimate $\hat{\mathbf B}=(\hat\beta_1,\ldots,\hat\beta_q)$ as $ {\mathbf B}^{(1)} $,
and the smoothing parameter $h_n$
used in obtaining it as $h_n^{(1)}$.
Construct a refined estimate of
$\nabla Q_\tau(\mathbf{x})$ as
\[
\hat\nabla Q^{(2)}_\tau(\mathbf{x}) =\bigl(\hat
c^{(2)}_{n,\u}(\mathbf {x};\tau )\bigr)_{\u\in A,[\u]=1} /
h_n^{(2)},
\]
where
%
%e3.7 #&#
\begin{equation}\quad
\hat c^{(2)}_{n,\u}(\mathbf{x};\tau) =\arg\min
_{\c} \sum_{i=1}^{n}
\rho_{\tau}\bigl\{Y_{i}-\mathbf{c}^\top\X
_{i\mathbf{x}}\bigl(h^{(1)}_{n},A\bigr)\bigr
\}K_{h^{(2)}_{n}}\bigl(\bigl|\X_{i\mathbf
{x}}^\top{\mathbf
B}^{(1)} \bigr|\bigr),\label{tech}
\end{equation}
and $K(\cdot)$ is a kernel density in $R^q$. Accordingly, the
estimates
$\hat\Sigma(\tau)$ and $\hat\Sigma_{\mathrm{T}}$ in~(\ref{jon}) and
(\ref{a6}) could be refined, respectively, as
\[
\hat\Sigma^{(2)}(\tau)= \frac{1}{n}\sum
_{j=1}^n \hat\nabla Q^{(2)}_\tau(
\X_j)\bigl\{\hat\nabla Q^{(2)}_\tau(
\X_j)\bigr\} ^\top
\]
and
\[
\hat\Sigma_{\mathrm{T}}^{(2)} = \int_{\delta^*}^{1-\delta^*}
\hat w^{(2)}(\tau) \hat\Sigma^{(2)}(\tau) \,d\tau,
\]
where $ \hat w^{(2)}(\tau) $
is constructed in the same way as $ \hat w(\tau) $, using eigenvalues
of $ \hat\Sigma^{(2)}(\tau) $.
Again, pick the first $ q $ eigenvectors of $ \hat\Sigma_{\mathrm{T}}^{(2)} $ to construct a new matrix ${\mathbf B}^{(2)}$ which can then
be substituted into (\ref{tech}) for ${\mathbf B}^{(1)}$. Repeat the
above two steps until convergence is reached. Intuitively, this
refined estimate of $\Sigma$ is more efficient due to the use of
a lower-dimensional kernel when estimating $\nabla Q_\tau(\mathbf{x})$,
thus mitigating the so-called ``curse of dimensionality'' problem.
We call the above procedure the adaptive composite quantile outer
product of gradients (qOPG).

We can also incorporate this ``composite-quantile'' idea
into the minimum average variance estimation (MAVE) procedure of \citet{r34} and propose a composite quantile MAVE (qMAVE) as
follows. With structural dimension $ q $, consider the following
minimization problem:
%
%e3.8 #&#
\begin{equation}
\int_{\delta^*}^{1-\delta^*} \sum_{j=1}^n
\sum_{i=1}^n \rho_{\tau}\bigl
\{Y_{i}- a_j - b_j^\top{\mathbf
B}^\top\X_{ij}\bigr\} K_{h_n}\bigl(|
\X_{ij}|\bigr) \,d \tau, \label{e15}
\end{equation}
with\vspace*{2pt} respect to $p\times q$ matrix ${\mathbf B}$, where $ \X_{ij} = \X_i
- \X_j $. Again, just as in (\ref{tech}), a possibly lower-dimensional kernel $K_{h_n}(|{\mathbf B}^\top\X_{ij}|)$ could be used to
replace $K_{h_n}(|\X_{ij}|)$ in (\ref{e15}), in the hope of an
improved efficiency of the resulted estimator, at least with
finite-sample size.
Estimates of the $q$ CS directions are thus given by the
orthonormalized columns of $ \hat{\mathbf B} $, the solution to
(\ref{e15}). Realization of (\ref{e15}) is similar to that of \citet{r35} and its theoretical properties can also be similarly
investigated by combining the results obtained in the \hyperref[app]{Appendix} and
the proofs in \citet{r35}.

To find out whether a qMAVE procedure would benefit from some
``adaptive'' weighting scheme,
one could consider, for example, a level-specific qMAVE procedure, where
\[
\hat{\mathbf B}(\tau)=\arg\min_{{\mathbf B}\in R^{p\times q}}\min_{a_j,b_j}
\sum_{j=1}^n \sum
_{i=1}^n \rho_{\tau}\bigl
\{Y_{i}- a_j - b_j^\top{\mathbf
B}^\top\X_{ij}\bigr\} K_{h_n}\bigl(|
\X_{ij}|\bigr) \,d \tau,
\]
and consequently define
\[
\hat\Sigma^*_{w} = \int_{\delta^*}^{1-\delta^*}
\hat{\mathbf B}(\tau )\hat{\mathbf B}(\tau)^\top\hat{w}(\tau) \,d\tau,
\]
where $\hat{w}(\tau)$ is the same as in (\ref{a6A}) derived from
the level-specific OPG matrix.
Our experience is such that this level-specific qMAVE is always
outperformed by both the qMAVE procedure of (\ref{e15})
and qOPG. A possible explanation is that $ \hat{\mathbf B}(\tau)$ being
an orthonormal matrix means that all directions [columns of $ \hat{\mathbf
B}(\tau)$] are equally weighted, whereas in qOPG the corresponding
directions (eigenvectors) are given different weights dictated by their
respective eigenvalues.

%s4 #&#
\section{Assumptions and theoretical results}\label{sec4}
For any $s_{0}=l+\gamma$, with
nonnegative integer $l$ and $0<\gamma\leq1$, we say a function
$m(\cdot)\dvtx R^{p}%
\rightarrow R$ has the order of smoothness $s_{0}$ on $\mathcal{D}$,
denoted by $m(\cdot)\in H_{s_{0}}(\mathcal{D})$ if, it is
differentiable up to order $l$ and there exists a constant $C>0$,
such that
\[
\bigl|D^{\mathbf{u}}m(\mathbf{x}_{1})-D^{\mathbf{u}}m(\mathbf
{x}_{2})\bigr|\leq C|\mathbf{x}_{1}-\mathbf{x}_{2}|^{\gamma}
\qquad\mbox{for all }\mathbf{x} %%
_{1},\mathbf{x}_{2}
\in\mathcal{D}\mbox{ and } {}[\mathbf{u}]=l.
\]

We assume the following conditions hold throughout the paper:
\begin{longlist}[(A1)]
\item[(A1)] The support $\mathcal{D}$ of $\X$ is open, convex and
the probability density function
of $\X$ is such that $f_{\X}(\cdot)\in H_{s_{1}}(\mathcal{D})$, for some
$s_1>0$.

\item[(A2)] The conditional quantile function
$Q_{\tau}(\cdot)\in H_{s_{2}}(\mathcal{D})$ for some $s_2>0$ uniformly
in $\tau\in(0,1)$.

\item[(A3)] There exist some positive values $\delta^*$, $ b_1$, $
b_2$ and $s_3>0$, such that
the conditional probability density $f_{Y|\X}(\cdot|\cdot)$ of $Y $ given
$\X$ belongs to $ H_{s_{3}}(\mathcal{D})$ and is uniformly bounded
away from zero in $(Q_\tau(\mathbf{x})-b_1,Q_\tau(\mathbf{x})+b_2)$
for all $\tau\in
[\delta^*,1-\delta^*]$ and $\mathbf{x}\in\mathcal{D}$.
%H_{s_4}(\mathcal{D})$, for some $s_4>0.$
\end{longlist}

The order of smoothness $s_1$, $s_2$, $s_3$ will be specified later.
The above assumptions are standard in local polynomial smoothing for
quantile regression; see, for example, \citet{r5} and
\citet{r19}. Among them, (A2) implies that for any $\mathbf
{x}\in
\mathcal{D} $ and $\X_i\in S_n(\mathbf{x})=\{i\dvtx 1\le i\le n,
|\X_{i\mathbf{x}}|\le h_n\}$, the error from approximating $Q_\tau
(\X_i)$
by the $k(=[s_{2}])$th order Taylor expansion
\[
Q_{n}(\X_i,\mathbf{x};\tau)=\bigl[\X_{i\mathbf{x}}(h_n,A)
\bigr]^\top\c _n(\mathbf{x};\tau)
\]
is of order $O(h_{n}^{s_{2}})$,
uniformly in $\{(\mathbf{x},\X_i)\dvtx
\mathbf{x}\in\mathcal{D}, \X_i\in S_n(\mathbf{x})\}$ and $\tau
\in(0,1)$.
(A3) strengthens Condition~3 of \citet{r5}, where it is
required that for a prespecified $\tau$,
$g(\mathbf{x}|\tau)=f_{Y|\X}(Q_{\tau}(\mathbf{x})|\mathbf{x})>0$,
for all
$\mathbf{x}\in\mathcal{D.}$

The following lemma concerns the strong uniform Bahadur type
representation of $\hat{\mathbf{c} }_{n}(\cdot;\tau)$ derived from
(\ref{a8}).

%
%le2 #&#
\begin{Lemma} \label{T1} Suppose \textup{(A1)--(A3)} hold with $s_1>0$, $s_2>0$,
$s_3>1/2$, and $k=[s_2]$.
The bandwidth $h_n$ is chosen such that
\[
h_n\propto n^{-\kappa}\qquad\mbox{with } \frac{1}{2(s_2+p)}\le
\kappa<\frac{1}{p}.
\]
Then we have with probability one
%
%e4.1 #&#
\begin{eqnarray}\label{correction}
&& \hat\c_{n}(\mathbf{x};\tau)-\c_{n}(\mathbf{x}; \tau)\nonumber
\\
&&\qquad =-\frac{\Sigma^{-1}_{n}(\mathbf
{x};\tau)}{N_{n}(\mathbf{x})}\sum_{i\in S_{n}(\mathbf{x})}
\X_{i\mathbf{x}}(h_{n},A)\bigl[I\bigl\{Y_{i}\le
Q_{n}(\X_{i},\mathbf{x};\tau)\bigr\}-\tau\bigr]
\\
&&\quad\qquad{}+O\biggl\{ \biggl(\frac{\log n}{nh_n^p} \biggr)^{ 3/4}\biggr\}\nonumber
\end{eqnarray}
uniformly in $\tau\in[\delta^*,1-\delta^*]$ and
$\mathbf{x}=\X_1,\ldots,\X_n$, where $N_n(\mathbf{x})=\sharp
S_n(\mathbf{x})$ and
\[
{\Sigma_{n}(\mathbf{x};\tau)}=\mathbf{E}_{i} \bigl[g(
\X%
_{i}|\tau)\X_{i\mathbf{x}}(h_{n},A)
\X_{i\mathbf{x}}^{\top
}(h_{n},A)|\X_i\in
S_n(\mathbf{x}) \bigr].
\]
\end{Lemma}
This strengthens the results obtained in \citet{r4} for
nonparametric quantile regression
and \citet{r19} for general nonparametric M-regression, both
of which concerned the uniformity in $\mathbf{x}$ only.
The uniformity in both $\mathbf{x}$ and $\tau$
plays a central role in examining the
asymptotic properties of $\hat\Sigma_{\mathrm{T}}$, defined via averaging over
$\mathbf{x}=\X_1,\ldots,\X_n$, and then integration with respect to
$\tau$ over $ [\delta^*,1-\delta^*]$.

We now move on to present the asymptotic properties of
$\hat\Sigma_{\mathrm{T}}$ and those of its eigenvalues and
eigenvectors. Write
$\nabla^2 Q_\tau(\cdot)$ for the Hessian matrix of $Q_\tau(\cdot)$ and
$\nabla g(\cdot|\tau)$, for the first-order derivative vector of $g(\cdot|\tau)$.
For any
$\tau\in(0,1)$ and $1\le k,l \le p$, let $\nabla
Q^{[k]}_{\tau}(\X)$ stand for the $k$th element of $\nabla
Q_{\tau}(\X)$; $\nabla^{[l]} g(\X|\tau)$ for
the $l$th element of $\nabla g(\X|\tau)$, $\nabla^2_{[k,l]}
Q_{\tau}(\X)$ for the $(k,l)$ element of $\nabla^2 Q_\tau(\cdot)$ and
write
\begin{eqnarray*}
&&\rho(\X|\tau,k,l)
\\
&&\qquad = \biggl[\frac{2\nabla^2_{[k,l]}
Q_{\tau}(\X)}{g(\X|\tau)}-\frac{\nabla Q^{[k]}_{\tau}(\X)
\nabla^{[l]} g(\X|\tau)}{g^2(\X|\tau)} -
\frac{ \nabla^{[l]}
Q_{\tau}(\X)\nabla^{[k]} g(\X|\tau)}{g^2(\X|\tau)} \biggr].
\end{eqnarray*}
For any $\tau_1,\tau_2\in(0,1)$ and $1\le k_1,l_1,k_2,l_2\le p$,
define
\begin{eqnarray*}
&& h(\tau_1,\tau_2|k_1,l_1,k_2,l_2)
\\
&&\qquad =\Cov \bigl( \nabla Q^{(k_1)}_{\tau_1}(\X)\nabla
Q^{(l_1)}_{\tau_1}(\X),\nabla Q^{(k_2)}_{\tau_2}(\X)
\nabla Q^{(l_2)}_{\tau_2}(\X) \bigr)
\\
&&\quad\qquad{}+\bigl\{\min(\tau_1,\tau_2)-\tau_1
\tau_2\bigr\} \Cov \bigl( \rho(\X|\tau,k_1,l_1),
\rho(\X|\tau,k_2,l_2) \bigr).
\end{eqnarray*}
For any symmetric $p\times p $ matrix ${\mathcal S}=(s_{ij})$, form a
$p(p+1)/2\times1$ vector
using the elements of ${\mathcal S}$:
\[
\operatorname{Vech}({\mathcal S})=(s_{11},\ldots,s_{p1},s_{22},
\ldots,s_{2p},s_{22},\ldots,s_{pp})^\top.
\]
Denote by $v(\cdot)$ the following $1 $-to-$1$ mapping from
$\{1,2,\ldots,p(p+1)/2\}$ onto $\{(i,j)\dvtx  1\le i\le j\le p\}$:
\[
v(k)=\bigl(v(k,1),v(k,2)\bigr)=(i,j)\qquad \mbox{such that }
\frac{(2p-i)(i-1)}{2}+j=k.
\]
In other words, the $k$th element of $\operatorname{Vech}({\mathcal S})$ is
given by
$s_{v(k)}=s_{v(k,1),v(k,2)}$.

Finally, for any symmetric $p\times p$ matrix ${\mathcal S}$, denote by
$\lambda_k({\mathcal S})$ and $\beta_k({\mathcal S})$, $k=1,\ldots,q$, the first $q$
(nonzero) eigenvalues and eigenvectors of ${\mathcal S}$, respectively. Write
$\tilde\lambda_{p-q}({\mathcal S})$ for the average of the smallest $p-q$
eigenvalues of ${\mathcal S}$.

%
%th1 #&#
\begin{Theorem} \label{T2}
Suppose \textup{(A1)--(A3)} hold with $s_1>0$, $s_3>1/2$, $s_2>3/2p+3$, and
$k=[s_2]$. Furthermore, the smoothing parameter $h_n$ is chosen such
that
%
%e4.2 #&#
\begin{equation}
h_n\propto n^{-\kappa}\qquad\mbox{with }\frac{1}{2(s_2-1)}\le
\kappa<\frac{1}{3p+4}. \label{mon}
\end{equation}
Then we have $\tilde\lambda_{p-q}(\hat\Sigma_{\mathrm{T}})=o_p(n^{-1/2})$ and
%
%e4.3 #&#
\begin{equation}
\sqrt{n}(\hat\Sigma_{\mathrm{T}}-\Sigma_{\mathrm{T}})\stackrel{d} {\to}
\mathbb{N},\label{donsker}
\end{equation}
where ``$\stackrel{d}{\to}$'' stands for convergence in distribution
and $\mathbb{N}$ stands for a
symmetric $p\times p$ random matrix, such that
$\operatorname{Vech}(\mathbb{N})$ is multivariate normal with zero mean and
covariance matrix ${\mathbf H}$,
whose $(k,l)$th element is given by
\[
\int_{0}^1\int_{0}^1
h\bigl(\tau_1,\tau_2|v(k,1),v(k,2),v(l,1),v(l,2)\bigr)
\,d\tau_1\,d\tau_2.
\]
Furthermore, if $\lambda_k(\Sigma_{\mathrm{T}})$, $k=1,\ldots,q$,
are all distinct, then for each $k=1,\ldots,q$,
%
%e4.4 #&#
%e4.5 #&#
\begin{eqnarray}
\sqrt{n}\bigl\{\lambda_k(\hat\Sigma_{\mathrm{T}})-
\lambda_k(\Sigma_{\mathrm{T}})\bigr\}&\stackrel{d} {\to}&
\beta^\top_k(\Sigma_{\mathrm{T}}) \mathbb{N}
\beta_k(\Sigma_{\mathrm{T}}), \label{evalue}
\\
\sqrt{n}\bigl\{\beta_k(\hat
\Sigma_{\mathrm{T}})-\beta_k(\Sigma_{\mathrm{T}})\bigr\}
&\stackrel{d} {\to}& \sum_{l=1,l\ne k}^q
\frac{\beta_l(\Sigma_{\mathrm{T}})\beta^\top
_{l}(\Sigma_{\mathrm{T}})
\mathbb{N}\beta_k(\Sigma_{\mathrm{T}})}{
\lambda_k(\Sigma_{\mathrm{T}})-\lambda_l(\Sigma_{\mathrm{T}})}. \label{evector}
\end{eqnarray}
\end{Theorem}
In theory, (\ref{evalue}) could be applied to make inference on the
structural dimension~$q$.
The proof of Theorem \ref{T2} is mainly based upon results on
U-processes [\citet{r27}], namely a collection of
U-statistics indexed by a family of symmetric kernels.

%s5 #&#
\section{Bandwidth selection}\label{sec5}\label{bandwidthsection}

As far as the point-wise estimation of $\nabla Q_\tau(\cdot)$ is
concerned, it followed from Lemma \ref{T1} that the ``optimal''
bandwidth $h_n$ which minimizes the pointwise\vspace*{1pt} mean square error
(MSE) of $\hat\nabla Q_\tau(\mathbf{x})$, is of the order $O(
n^{-1/(p+2k+2)})$. In this sense, the choice (\ref{mon}) of the
bandwidth $h_n$
under-smooths the estimator. Such undersmoothing is necessary for the
estimator $\hat\nabla Q_\tau(\mathbf{x})$ to have a bias of order $o_p(
n^{-1/2})$ thus negligible. The stochastic term of $\hat\nabla
Q_\tau(\mathbf{x})$, once averaged over $\mathbf{x}=\X_1,\ldots,\X
_n$, can achieve
the rate of $O_p(n^{-1/2})$, independent of the speed at which $h_n$
tends to zero. Similar observations have been made in \citet{r5} and \citet{r20}. In cases where the link function
$Q_{\tau}(\cdot)$ closely resembles a (local)
polynomials,\vadjust{\goodbreak}
the bias thus becomes less of an issue as it either significantly
reduces or completely vanishes;
we can then afford to employ a larger bandwidth, and thus produce more
efficient estimates of $\nabla Q_\tau(\cdot)$,
while results in Theorem \ref{T2} still hold. This also
explains our assertion in Section~\ref{sec1} that qOPG is
structure-adaptive.
In practice, an empirical ``optimal'' bandwidth can be obtained by
plugging in estimates for the unknown quantities in the formula of
the pointwise theoretical ``optimal'' bandwidth.

We can also select bandwidth based on the cross-validation (CV)
criterion for quantile regression %;see, for example, Al-kenani and Yu (\citeyear{AY2010}). This is carried out
as follows.
For any given $ \tau\in(0,1) $ and fixed $h_n$,
denote by $ Q^{\setminus j}_\tau(\mathbf{x}|h_n)$, $ j=1,\ldots,n$,
the leave-one-out estimate of $ Q_\tau(\X_j)$
using $ \{(X_i, Y_i)\dvtx  i \neq j\} $ with bandwidth $h_n$. Let
\[
\mathrm{CV}(\tau, h_n) = n^{-1} \sum_{j=1}^n
\rho_\tau \bigl(Y_j - Q^{\setminus j}_\tau(X_j|h_n)
\bigr),
\]
and denote by $ h_\tau^{\mathrm{CV}} $, the level-specific cross-validated (CV)
bandwidth,
namely the $ h_n$ that minimizes $ \mathrm{CV}(\tau, h_n) $.
However, based on our experience with simulated data,
we found such level-specific CV bandwidth selection is not only rather
time-consuming,
but also terribly unstable,
possibly due to the difficulty in assessing the goodness-of-fit
in quantile regression; see \citet{r15}.
Instead, we recommend the following modified level-specific CV bandwidth.
First, consider an average of the level-specific CV bandwidth $
h_\tau^{\mathrm{CV}} $ with $ \tau$ ranging over the set of
$\{\tau_s = s/(T+1)\dvtx  s = 1, \ldots, T\} $ for some positive integer $T$:
\[
\bar h^{\mathrm{CV}} = \sum_{s=1}^T
h_{\tau_s}^{\mathrm{CV}} /T.
\]
Then in view of the relationship proposed in \citet{r39}, we
define the modified level-specific CV bandwidth as
%
%e5.1 #&#
\begin{equation}
\bar h_\tau^{\mathrm{CV}} = \bar h^{\mathrm{CV}} \bigl\{ \tau(1-
\tau)/\phi\bigl(\Phi^{-1}(\tau )\bigr)\bigr\}^{1/5},\label{cv}
\end{equation}
where functions
$\phi(\cdot) $ and $ \Phi(\cdot) $ are, respectively, the probability
and cumulative distribution functions of the standard normal distribution.
Compared to $ h_\tau^{\mathrm{CV}} $, $\bar h_\tau^{\mathrm{CV}}$ is more stable and
delivers much better results,
but its computation is equally\vadjust{\goodbreak} computationally
intensive. We also tried out variations of $\bar h_\tau^{\mathrm{CV}}$
defined as in (\ref{cv}) but with
$ \bar h^{\mathrm{CV}} $ replaced by bandwidths chosen via other procedures.
Our best experience lies with $\bar h_\tau^{\mathrm{CV}}$ with $\bar h^{\mathrm{CV}}$
set to be the CV bandwidth for conditional mean
regression of $ |Y - E(Y)| $ on $ \X$.

%s6 #&#
\section{Estimation of the structural dimension}\label{sec6}\label{qsection}
According to Theorem \ref{T2}, the average of the smallest $p-q$
eigenvalues of $\hat\Sigma_{\mathrm{T}}$ defined in (\ref{a6}) is of
order $o_p(n^{-1/2})$. For $k=1,\ldots,p$, plot the average of the
smallest $k$ eigenvalues of $\hat\Sigma_{\mathrm{T}}$ against $k$ and
likely values for $q$ could be then identified by noting the
location of a noticeable increase.
The asymptotic distribution of the eigenvalues of $ \hat\Sigma_{\mathrm{T}} $
given in Theorem \ref{T2} could also be used for selecting $ q$. However,
as the distribution depends on another unknown matrix $ {{\mathbf H}} $
which is not easy to estimate, such approach might not be very
practical.

Combining the CV method of \citet{r34} with the composite
quantile regression provides an alternative way to select $ q $. For
illustration purposes, we here give details for the local constant
quantile kernel smoothing.
With working dimension $ q $, suppose the $q$-columns of $ \hat B_q $
are the corresponding estimates of the~CS directions.
For each observation $ (X_j, Y_j), j = 1, \ldots, n $,
calculate the delete-one-estimator of $ \tilde Q_\tau(\hat B_q X_j ) $ of
(\ref{yuanquan}) as
\[
\hat Q^{\setminus j}_\tau\bigl(\hat B_q^\top
X_j\bigr) = \arg\min_c \sum
_{i\neq j} \rho_\tau(Y_i - c)
K_{h_n}\bigl(|\hat B_q X_{ij}|\bigr).
\]
We then define the CV value specific to working dimension $ q $ as
\[
\mathrm{CV}(q) = \int_{\delta^*}^{1-\delta^*}\sum
_{j=1}^n \rho_\tau\bigl(Y_i
- \hat Q^{\setminus j}_\tau\bigl(\hat B_q^\top
X_j\bigr)\bigr) \,d\tau,
\]
and choose the dimension which minimizes $ \mathrm{CV}(q) $. Our simulation
study suggests that this methodology works reasonably well, though
it is also rather computationally intensive.

%s7 #&#
\section{Numerical study}\label{sec7}

In this section, we first carry out comparison studies of
the two newly proposed procedures, qOPG and qMAVE, with two existing methods
using
simulated data. The two new procedures are then applied to the analysis
of a real data set for the purpose of discovering the dimension
reduction space.

%t1 #&#
\begin{table}
\def\arraystretch{0.99}
\tabcolsep=0pt
\caption{Average estimation errors and their standard derivation (in parenthesis)
and frequency of correct structural dimension identification}\label{tab1}
\begin{tabular*}{\tablewidth}{@{\extracolsep{\fill}}@{}lcc cccccccc @{}}
\hline
& & & & & &\multicolumn{2}{c}{\textbf{qOPG}} & \multicolumn{2}{c}{\textbf{qMAVE}}\\[-6pt]
& & & & & &\multicolumn{2}{c}{\hrulefill} & \multicolumn{2}{c}{\hrulefill}\\
\textbf{Model} & $\bolds{\varepsilon}$ & $\bolds{n}$ & \textbf{SIR} & \textbf{dOPG} & \textbf{dMAVE} & $\bolds{h_0}$ & $\bolds{h_{\mathrm{CV}}}$ & $\bolds{h_0}$
& $\bolds{h_{\mathrm{CV}}}$ & \textbf{freq.}\\
\hline
(A) & $N(0,1)$ & 200 & 0.82 & 0.55 & 0.53 & 0.42 & 0.44 & 0.48 & 0.48 &56\%\\
& &&{(0.14)}& {(0.20)} & {(0.18)} & {(0.15)} & {(0.15)} & {(0.16)} & {(0.15)} &\\
& & 400 &0.68& 0.37 & 0.35 & 0.27 & 0.26 & 0.31 & 0.30 & 90\%\\
& &&{(0.16)} & {(0.14)} & {(0.10)} & {(0.08)} & {(0.08)} & {(0.08)} & {(0.08)}&\\
& $t(3)/\sqrt{3}$ & 200 &0.79& 0.50 & 0.46 & 0.42 & 0.38 & 0.38 & 0.40&72\%\\
& &&{(0.15)} & {(0.22)} & {(0.16)} & {(0.15)} & {(0.14)} & {(0.14)} & {(0.14)}&\\
& & 400 &0.63& 0.31 & 0.29 & 0.22 & 0.21 & 0.23 & 0.24 & 97\%\\
& &&{(0.16)} & {(0.13)} & {(0.08)} & {(0.07)} & {(0.07)} & {(0.06)} & {(0.06)}&\\
& $\chi^2(1)$ & 200 & 0.78& 0.61 & 0.50 & 0.48 & 0.49 & 0.46 & 0.49& 50\%\\
& &&{(0.13)} & {(0.22)} & {(0.17)} & {(0.20)} & {(0.19)} & {(0.17)} & {(0.17)}&\\
& & 400 &0.61& 0.39 & 0.32 & 0.30 & 0.28 & 0.28 & 0.29 & 79\%\\
& &&{(0.14)} & {(0.16)} & {(0.10)} & {(0.12)} & {(0.09)} & {(0.09)} & {(0.10)}&
\\[3pt]
(B) & $N(0,1)$ & 200 &0.69& 0.58 & 0.59 & 0.44 & 0.50 & 0.54 & 0.52 &56\%\\
& &&{(0.17)} & {(0.17)} & {(0.18)} & {(0.18)} & {(0.19)} & {(0.19)} & {(0.19)}&\\
& & 400 &0.51& 0.35 & 0.38 & 0.24 & 0.27 & 0.32 & 0.32 & 87\%\\
& &&{(0.15)} & {(0.10)} & {(0.13)} & {(0.10)} & {(0.10)} & {(0.11)} & {(0.12)}&\\
& $t(3)/\sqrt{3}$ & 200 &0.57& 0.48 & 0.47 & 0.38 & 0.37 & 0.40 & 0.40 &84\%\\
& &&{(0.16)} & {(0.16)} & {(0.15)} & {(0.16)} & {(0.12)} & {(0.15)} & {(0.13)}&\\
& & 400 &0.41& 0.34 & 0.29 & 0.19 & 0.18 & 0.21 & 0.22 & 97\%\\
& &&{(0.12)} & {(0.10)} & {(0.09)} & {(0.09)} & {(0.06)} & {(0.06)} & {(0.07)}&\\
& $\chi^2(1)$ & 200 &0.64& 0.57 & 0.53 & 0.55 & 0.46 & 0.51 & 0.48 &64\%\\
& &&{(0.17)} & {(0.18)} & {(0.20)} & {(0.24)} & {(0.20)} & {(0.22)} & {(0.19)}&\\
& & 400 &0.42& 0.35 & 0.31 & 0.24 & 0.22 & 0.24 & 0.25 & 94\%\\
& &&{(0.11)} & {(0.13)} & {(0.09)} & {(0.11)} & {(0.08)} & {(0.07)} & {(0.07)}&
\\[3pt]
(C) & $N(0,1)$ & 200 &0.53& 0.55 & 0.51 & 0.77 & 0.42 & 0.48 & 0.36 &29\%\\
& &&{(0.13)} & {(0.14)} & {(0.17)} & {(0.15)} & {(0.14)} & {(0.17)} & {(0.10)}&\\
& & 400 &0.37& 0.36 & 0.33 & 0.77 & 0.29 & 0.30 & 0.24 & 31\%\\
& &&{(0.08)} & {(0.11)} & {(0.09)} & {(0.16)} & {(0.10)} & {(0.08)} & {(0.05)}&\\
& $t(3)/\sqrt{3}$ & 200 &0.61& 0.62 & 0.59 & 0.81 & 0.47 & 0.55 & 0.38 &32\%\\
& &&{(0.15)} & {(0.14)} & {(0.18)} & {(0.15)} & {(0.19)} & {(0.19)} & {(0.14)}&\\
& & 400 &0.44& 0.41 & 0.38 & 0.77 & 0.39 & 0.35 & 0.25 & 39\%\\
& &&{(0.12)} & {(0.14)} & {(0.15)} & {(0.15)} & {(0.20)} & {(0.14)} & {(0.07)}&\\
& $\chi^2(1)$ & 200 &0.63& 0.60 & 0.49 & 0.50 & 0.46 & 0.44 & 0.42 &37\%\\
& &&{(0.14)} & {(0.15)} & {(0.16)} & {(0.17)} & {(0.16)} & {(0.16)} & {(0.13)}&\\
& & 400 &0.43& 0.42 & 0.32 & 0.35 & 0.30 & 0.31 & 0.27 & 46\%\\
& &&{(0.11)} & {(0.14)} & {(0.09)} & {(0.10)} & {(0.16)} & {(0.09)} & {(0.08)}&\\
\hline
\end{tabular*}\vspace*{-3pt}
\end{table}

In the calculation below, the local linear quantile regression, that
is, $ k = 1 $, and the Epanechnikov kernel
function are used. The integrations in (\ref{a6A}) and (\ref{e15})
are evaluated by the weighted summation
of $ \hat\Sigma(\tau) $ over $ \tau= 0.1, 0.2, \ldots, 0.9 $.

%
%ex1 #&#
\begin{Example}[(Simulated data)] \label{Example1}
We reconsider the following three models that are commonly tested
out in the field of dimension reduction:
\begin{eqnarray*}
\mbox{Model (A):}\qquad Y &=& \mathbf{x}_1(\mathbf{x}_1+
\mathbf{x}_2+1) + 0.5\varepsilon,
\\
\mbox{Model (B):}\qquad Y &=& \mathbf{x}_1/\bigl(0.5 + (
\mathbf{x}_2+1.5)^2\bigr) + 0.5\varepsilon,
\\
\mbox{Model (C):}\qquad Y &=& \mathbf{x}_1 + \exp(
\mathbf{x}_2) \varepsilon,
\end{eqnarray*}
where $ \X= (\mathbf{x}_1, \ldots, \mathbf{x}_{10})^\top\sim N(0,
(\sigma_{ij})_{1\le i,j \le10} ) $ with $ \sigma_{ij} = 0.5^{|i-j|}
$, and $\varepsilon$ is the error term designed to have various
distributions; see Table~\ref{tab1} below.
The first two models were thoughtfully designed by \citet{r23}
for the study of Slice Inverse Regression (SIR).
Model 3 was used in \citet{r35} in the context of conditional mean
and conditional variance based dimension reduction.

Based on the conclusion of \citet{r25} from their intensive
comparison study using simulated data, we have chosen to compare our
conditional quantile-based approaches, qOPG and qMAVE, with dOPG
and dMAVE of \citet{r35}, among the many existing dimension
reduction procedures. Another reason for us to include dOPG and
dMAVE in the study is the fact that these conditional
probability-based approaches, are theoretical equivalences to qOPG
and qMAVE. We hope through such comparison can manifest the
structure-adaptive nature of our new methods. We also include in the
comparison study the SIR of \citet{r23}, for which 8 slices are used
when the sample size $n=200$, and 10 when the sample size $n=400$.
For dOPG and dMAVE, following the \mbox{rule-of-thumb} as in \citet{r35},
we use bandwidths of order $n^{-1/5} $ and $n^{-1/(p+4)}$,
respectively, for the two kernels in the estimation. For qOPG and
dMAVE, the bandwidth is chosen as described in Section~\ref{bandwidthsection}. For any estimator $ \hat{\mathbf B} $ of $ {\mathbf
B}_0 $, we define the estimation error as the largest among the
absolute values of the elements of $
\hat{\mathbf B} (\hat{\mathbf B}^\top\hat{\mathbf B} )^{-1} \hat{\mathbf B} -
{\mathbf B}_0 ( {\mathbf B}_0^\top
{\mathbf B}_0 )^{-1} {\mathbf B}_0$. Table~\ref{tab1} reports the mean and standard error (in brackets) of the
estimation error from 100 replicates for various combinations of
model, error distribution and sample size.
The last column of Table~\ref{tab1}
is the percentage of times that the structural dimension has been
correctly identified by the CV method described in Section~\ref{qsection}.

A general observation is such that qOPG and qMAVE---either with
data-driven bandwidth or with
a bandwidth chosen according to the rule-of-thumb---outperform,
respectively, dOPG and dMAVE as well as SIR for both models (A)~and~(B). The only exception lies with model (C), where qOPG using the
rule-of-thumb bandwidth is beaten by dOPG, but the situation
reverses with a data-driven bandwidth.
This provides a line of empirical evidence for the assertion we made in
Section~\ref{sec1} that
if the conditional quantile function is well approximated locally by
polynomials, then the data-driven bandwidth deduced from qOPG means
more efficient estimators. Another noticeable pattern is that,
contradictory to what happens with conditional density-based methods
where dMAVE consistently outperforms dOPG, the expected superiority
of qMAVE over qOPG is nowhere obvious. In fact, for models (A) and
(B), qOPG outperforms qMAVE most of the time, especially so when
data-driven bandwidths are used. Even for model (C) qMAVE seems to
enjoy an obvious lead over qOPG, this again becomes less obvious
when a data-driven bandwidth is used. A plausible explanation for this
might be that an
adaptive-weighting scheme has been incorporated into qOPG,
while such procedure is hard to be combined with qMAVE.
\end{Example}

%
%ex2 #&#
\begin{Example}[(Real data)] \label{Example2}
In financial economics, the capital asset pricing model (CAPM)
indicates that the return
of a portfolio strongly depends on the market performance.
However, little is known about the factors that affect the\vadjust{\goodbreak} volatility
of a portfolio.
In the following, we consider the daily return $ Y $ of a portfolio
listed~at
\fontsize{9.5pt}{\baselineskip}\selectfont
\[
\mbox{\url{http://mba.tuck.dartmouth.edu/pages/faculty/ken.french/data\_library.html}}
\]
\normalsize
with covariate $ \X= (\mathbf{x}_1, \mathbf{x}_2, \ldots, \mathbf
{x}_{15}) $, where $ \mathbf{x}_1, \ldots, \mathbf{x}_5 $
are the returns of the portfolio in the past five days, and
$ \mathbf{x}_6, \ldots, \mathbf{x}_{10} $ are the absolute values of the
returns which are proxy of the past volatilities;
$ \mathbf{x}_{11}, \ldots, \mathbf{x}_{15} $ are the market returns on
the same day as $ Y $ and those in the past four days,
and $ \mathbf{x}_{16}, \ldots, \mathbf{x}_{20} $ are the absolute values
of the market returns.

%
%t2 #&#
\begin{table}
\tabcolsep=0pt
\caption{Estimated CS directions for Example \protect\ref{Example2}}\label{tab2}
\begin{tabular*}{\tablewidth}{@{\extracolsep{\fill}}@{}lcc ccc ccc ccc@{}}
\hline
$\bolds{\mathbf{x}_i}$ & $\bolds{\beta_1}$ & $\bolds{\beta_2}$ &
$\bolds{\mathbf{x}_i}$ & $\bolds{\beta_1}$ & $\bolds{\beta_2}$ &
$\bolds{\mathbf{x}_i}$ & $\bolds{\beta_1}$ & $\bolds{\beta_2}$ &
$\bolds{\mathbf{x}_i}$ & $\bolds{\beta_1}$ & $\bolds{\beta_2}$\\
\hline
$\mathbf{x}_1 $ & $-$0.014 & \phantom{$-$}0.001 & $\mathbf{x}_6$ & \phantom{$-$}0.006 & $-$0.089 & $\mathbf{x}_{11} $ & \phantom{$-$}0.994 & $-$0.032 & $\mathbf{x}_{16} $ & \phantom{$-$}0.029 & \phantom{$-$}0.490\\
$\mathbf{x}_2 $ & $-$0.042 & \phantom{$-$}0.045 & $\mathbf{x}_7 $ &\phantom{$-$}0.005 & \phantom{$-$}0.093 & $\mathbf{x}_{12} $ & \phantom{$-$}0.048 & $-$0.027 & $\mathbf{x}_{17} $ & $-$0.017 & \phantom{$-$}0.506\\
$\mathbf{x}_3 $ & $-$0.029 & $-$0.239 & $\mathbf{x}_8 $ &\phantom{$-$}0.020 & \phantom{$-$}0.271 & $\mathbf{x}_{13} $ & \phantom{$-$}0.048 & $-$0.076 & $\mathbf{x}_{18} $ & \phantom{$-$}0.008 & \phantom{$-$}0.302\\
$\mathbf{x}_4 $ & $-$0.008 & $-$0.100 & $\mathbf{x}_9 $ &\phantom{$-$}0.008 & \phantom{$-$}0.277 & $\mathbf{x}_{14} $ & \phantom{$-$}0.035 & \phantom{$-$}0.347 & $\mathbf{x}_{19} $ & $-$0.035 & $-$0.126\\
$\mathbf{x}_5 $ & \phantom{$-$}0.008 & \phantom{$-$}0.067 & $\mathbf{x}_{10} $ &$-$0.014 & \phantom{$-$}0.111 & $\mathbf{x}_{15} $ & $-$0.005 & \phantom{$-$}0.010 & $\mathbf{x}_{20} $ & \phantom{$-$}0.001 & $-$0.120\\
\hline
\end{tabular*}
\end{table}

%
%
%f1 #&#
\begin{figure}%[hp!]

\includegraphics{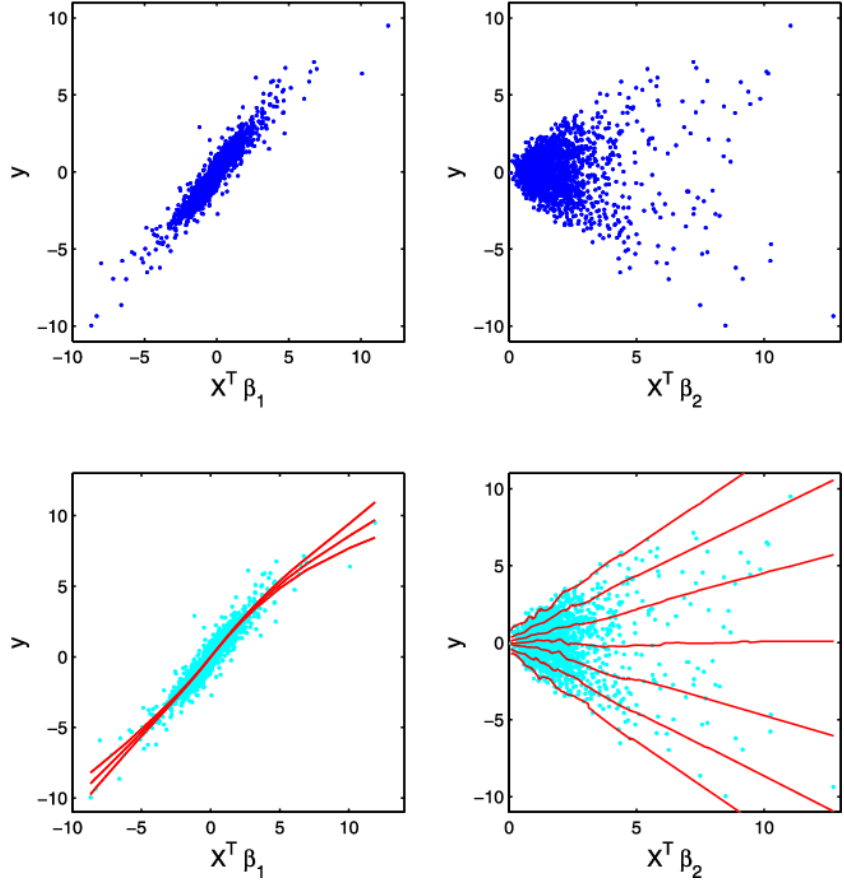}

\caption{Results for Example \protect\ref{Example2}. The top
two panels are the scatter plots of $ Y $ against the two estimated
CS directions $\beta_1$ and $\beta_2$. The bottom-left panel is the
fitted regression function of $ Y $ against the first CS direction
and its 95\% confidence interval.
In the bottom-right panel, the curves are the regression quantiles of $
Y $ against the second directions at
$ \tau= 0.01, 0.1, 0.3, 0.5, 0.7, 0.9, 0.99$, respectively.} \label{figExample}
\end{figure}

Applying qOPG, the first several eigenvalues of $ \hat\Sigma_{\mathrm{T}} $ are, respectively, $1.0620$, $0.0164$,
$0.0017$, $0.0007$ and $0.0004$. With the structural dimension
set as $ 2$, we
obtain the estimated CS directions $\beta_1$ and $\beta_2$; see Table~\ref{tab2}.
The scatter plots of $ Y $ against $\beta_1^\top\X$ and $\beta
_2^\top\X$
are given in Figure~\ref{figExample}. The fitted curve in the bottom
two panels are
created
with bandwidths $ h = h_0/(\hat f_k(x))^{0.2} $ with $ h_0 $ being
selected by the CV method,
and $ \hat f_k(\cdot)$, $k = 1, 2$, being the kernel estimate of the
density function of $ \beta_k^\top\X$.
The fitted regression function of the portfolio's return on $
\beta_1^\top\X$ in the bottom-left panel of Figure~\ref{figExample} suggests the first CS direction $\beta_1$ is mostly
about the conditional mean,
while the second CS direction $\beta_2$ is clearly about the
conditional variance, evident from
the bottom-right panel.
The first direction $ \beta_1 $ is dominated by $ \mathbf{x}_{11}$,
the market return of the day,
with a coefficient 0.9940; this is in line with the CAPM in that the
expected return of any portfolio
largely depends on the present-day market performance. It is also
interesting to note that the volatility of the portfolio also
depends the market's volatility,
as suggested by the large coefficients of $ \mathbf{x}_{16}, \mathbf
{x}_{17} $ and $ \mathbf{x}_{18} $
on the second CS direction $ \beta_2 $.
Also, its own past volatilities
($ \mathbf{x}_8, \mathbf{x}_9 $)
also contribute to its present-day volatility, although to a less extent.
\end{Example}

%s8 #&#
\section{Conclusions}\label{sec8}
In this paper, we have proposed and investigated two composite
quantile approaches to dimension reduction, namely qOPG and qMAVE.
Compared with moment-based methods, these methods require less
restrictive assumptions and can identify all dimension reduction
directions. It does not involve ``slicing'' of the response variable
$ Y$, as is the case with SIR or conditional density-based methods
[\citet{r35}].
It carries out regression analysis directly on $ Y $
instead of transformations of $ Y $.
As a result of these characteristics, qOPG and qMAVE are
structure-adaptive, and thus more efficient.
However, because the amount of computation embedded in quantile
regression is significantly heavier than in least
square minimization, the implementation of qOPG and qMAVE is rather
time consuming compared to
most of the existing methods. Because of this, we recommend the use of
dOPG or dMAVE to obtain an initial estimator of the central
subspace and of the structural dimension, and the use of qOPG or
qMAVE for more efficient refined estimator.

%sA #&#
\begin{appendix}\label{app}
\section*{Appendix: Proofs}
\renewcommand{\theLemma}{A.\arabic{Lemma}}
\setcounter{Lemma}{0}

\begin{pf*}{Proof of Lemma \ref{lemma1}}
The assertion that ${\mathcal S}(\Sigma)\subseteq{\mathcal S}({\mathbf
B}_0)$ follows directly from~(\ref{a2}).
We show that
the opposite holds too.
Based on (\ref{a2}), we can see by definition
\[
\Sigma={\mathbf B}_0 \biggl[\int_0^1
E\bigl\{\nabla\tilde Q_\tau\bigl({\mathbf B}_0^\top
\X\bigr) \bigl[\nabla\tilde Q_\tau\bigl({\mathbf B}_0^\top
\X\bigr) \bigr]^\top\bigr\} \,d\tau \biggr] {\mathbf
B}_0^\top.
\]
It thus suffices if we can prove the matrix
\[
M=\int_0^1 E\bigl\{\nabla\tilde
Q_\tau\bigl({\mathbf B}_0^\top\X\bigr)\bigl[
\nabla\tilde Q_\tau\bigl({\mathbf B}_0^\top\X
\bigr) \bigr]^\top\bigr\} \,d\tau
\]
is of full rank. For if otherwise, there must exist some vector
${\mathbf b}_1\in R^q$, with Euclidean norm one
such that ${\mathbf b}_1^\top M{\mathbf b}_1=0$. Seeing the definition
of $M$, this implies that
%
%eA.1 #&#
\begin{equation}
{\mathbf b}_1^\top\nabla\tilde Q_\tau\bigl({
\mathbf B}_0^\top\X\bigr)=0\qquad\mbox{a.s.} \label{a3}
\end{equation}
for all $\tau\in(0,1)$ except on a set of Lebesgue measure zero.

Let $\mathbf{B}=({\mathbf b}_1,\ldots,{\mathbf b}_q)\in R^{q\times
q}$ denote an
orthonormal basis for $R^q$, that is, $\mathbf{B}^\top\mathbf
{B}={\mathbf I}_q$.
For any given $\tau\in(0,1)$, write
%
%eA.2 #&#
\begin{equation}
G_\tau(\u) = \tilde Q_\tau(\u),\qquad\tilde
G_\tau(\u) = \tilde Q_\tau(\mathbf{B}\u), \qquad\tilde{
\mathbf B}_0={\mathbf B}_0\mathbf{B}.\label{a4}
\end{equation}
Thus,
\[
G_\tau({\mathbf B}\u)=\tilde G_\tau(\u);\qquad
G_\tau\bigl({\mathbf B}_0^\top\X \bigr) =
G_\tau\bigl(\tilde{\mathbf B}_0^\top\X\bigr).
\]
Consider the gradient
vector of $\tilde G_\tau(\u)$ and then evaluate it for $\u=\tilde
{\mathbf B}_0^\top\X$:
\begin{eqnarray*}
\frac{\partial\tilde G_\tau(\u)}{\partial\u}&=&\frac{\partial
G_\tau({\mathbf B}\u)}{\partial\u} ={\mathbf B}^\top
\frac{\partial G_\tau({\mathbf B}\u)}{\partial({\mathbf B}\u
)}={\mathbf B}^\top\nabla G_\tau({\mathbf B}
\u) \stackrel{\u=\tilde{\mathbf B}_0^\top\X}=
\mathbf{B}^\top\nabla G_\tau\bigl( {\mathbf
B}_0^\top\X\bigr),
\end{eqnarray*}
the first element of which, according to (\ref{a3}), equals zero.
This suggests the value of $\tilde G_\tau(\tilde{\mathbf B}_0^\top\X)$,
as a function of
$\tilde{\mathbf B}_0^\top\X=({\mathbf b}_1^\top{\mathbf B}_0^\top\X,\ldots, {\mathbf b}_q^\top{\mathbf B}_0^\top\X)^\top$,
does not change with ${\mathbf b}_1^\top{\mathbf B}_0^\top\X$.
This together with the fact that
\[
\tilde G_\tau\bigl(\tilde{\mathbf B}_0^\top\X
\bigr)=G_\tau\bigl({\mathbf B}_0^\top\X\bigr)=
\tilde Q_\tau\bigl({\mathbf B}_0^\top\X\bigr)=
Q_\tau( \X)
\]
implies that $ Q_\tau( \X)$ is in fact a function of $q-1$ variables:
${\mathbf b}_2^\top{\mathbf B}_0^\top\X,\ldots,\break  {\mathbf b}_q^\top{\mathbf
B}_0^\top\X$.
And this according to (\ref{a3}) holds for any $\tau\in(0,1)$.
As $\{Q_\tau( \X)\dvtx \tau\in
(0,1)\}$ collectively defines $F(\cdot|\X) $, we can conclude that $F(\cdot|\X
) $
is in fact a function of $({\mathbf b}_2^\top{\mathbf B}_0^\top\X,\ldots, {\mathbf b}_q^\top{\mathbf B}_0^\top\X)^\top
=[{\mathbf B}_0({\mathbf b}_2,\ldots,{\mathbf b}_q)]^\top\X$,
expressed as
\[
F(Y| \X) = F\bigl(Y|\tilde{\mathbf B}^\top\X\bigr),\qquad\mbox{a.s.
where }\tilde{\mathbf B}= {\mathbf B}_0({\mathbf b}_2,
\ldots,{\mathbf b}_q).
\]
This means $S(\tilde{\mathbf B})$ is SDR and as $S({\mathbf B}_0)$ is the CS,
we should have
$S({\mathbf B}_0)\subseteq S(\tilde{\mathbf B}) $. This contradicts the fact that
$\dim(S({\mathbf B}_0))=q>q-1=\dim( S(\tilde{\mathbf B})) $.
\end{pf*}

The proof of Lemma \ref{T1} is left until the end.
To prove Theorem \ref{T2}, we also need to introduce more notation.
For any
$\t=(t_1,\ldots,t_p)^\top\in[-1,1]^{p}$, let $\mathbf{t}(A)$ stand
for the $s(A)\times1$
vector $(\mathbf{t}^{\mathbf{u}})_{\mathbf{u}\in A}$.
Define
\[
\Gamma=\int_{[-1,1]^{p}}\t(A)\bigl\{\t(A)\bigr\}^\top \,d
\t.
\]
Standard result in kernel smoothing [e.g., \citet{r26}] is such that
with probability one,
%
%eA.3 #&#
\begin{equation}
\frac{N_{n}(\mathbf{x})}{nh_n^p}-f_{\X}(\mathbf {x})=O\bigl(h_n^2+
\bigl(nh^p_n/\log n\bigr)^{-1/2}\bigr)
\label{ellen1}
\end{equation}
uniformly in $ \mathbf{x}\in\mathcal{D}$, and
%
%eA.4 #&#
\begin{equation}
\Sigma_n(\mathbf{x};\tau)-g(\mathbf{x}|\tau)\Gamma =O \bigl(
\bigl(nh_n^p/\log n\bigr)^{-1/2}+h_n
\bigr) \label{ellen}
\end{equation}
uniformly in $\tau\in(0,1)$ and $\mathbf{x}\in\mathcal{D}$.

Also, we will cite the following result, the proof of which will be
given at the end of this section:
with probability one,
%
%eA.5 #&#
\begin{eqnarray}\label{lee}
\qquad&& \sum_{i}
\X_{i\mathbf{x}}\bigl(h_{n},A\bigr)I\bigl(|
\X_{i\mathbf{x}}|\le h_n\bigr) \bigl[I\bigl\{Y_{i}\le
Q_{n}(\X_{i},\mathbf{x};\tau)\bigr\}
-I\bigl\{Y_{i}\le Q_\tau(\X_{i})\bigr\}
\bigr]
\nonumber\\[-8pt]\\[-8pt]
&&\qquad =o\bigl(n^{-1/2}\bigr) \nonumber
\end{eqnarray}
uniformly in $\mathbf{x}\in\mathcal{D}$, $\tau\in(0,1)$.

\begin{pf*}{Proof of Theorem \ref{T2}}
Write as $\tilde\Gamma_n(\X_j;\tau)$, the $p\times s(A)$ matrix
consisting the second up to the $(p+1)$th row of $\Sigma^{-1}_{n}(\X
_j;\tau)$.
First note that under conditions in Theorem \ref{T2},
\begin{eqnarray*}
h_n^{-1}\bigl(nh_n^p/\log n
\bigr)^{-3/4}&=&o\bigl(n^{-1/2}\bigr),\qquad h_n^{s_2-1}=o
\bigl(n^{-1/2}\bigr)\quad\mbox{and}
\\
\log n/\bigl(nh_n^p\bigr)&=&o\bigl(n^{-1/2}h_n\bigr).
\end{eqnarray*}
This together with (\ref{jon}) and Lemma \ref{T1} leads to
\begin{eqnarray*}
\hat\Sigma(\tau)&=&\frac{1}{n}\sum_{j=1}^n
\nabla Q_\tau(\X_j)\bigl\{\nabla Q_\tau(
\X_j)\bigr\}^\top +h_n^{-1}
\bigl[M_n(\tau)+M_n^\top(\tau)\bigr] +o
\bigl(n^{-1/2}\bigr),
\end{eqnarray*}
where
\begin{eqnarray*}
M_n(\tau)&=&\frac{1}{n}\sum_{i,j}
\frac{\nabla Q_\tau(\X
_j)}{N_n(\X_j)} I\bigl(|\X_{ij}|\le h_n\bigr)\bigl[I\bigl
\{Y_{i}\le Q_{n}(\X_{i},\X_j;\tau)
\bigr\}-\tau\bigr]
\\
&&\hspace*{23pt}{}\times \X^\top_{ij}(h_{n},A)\tilde
\Gamma _{n}^\top(\X_j;\tau)
\end{eqnarray*}
with $\X_{ij}=\X_{i}-\X_j$. Using results in (\ref{ellen1}), (\ref
{ellen}) and (\ref{lee}),
we have
%
%eA.6 #&#
\begin{eqnarray}\label{ling}
\hat\Sigma(\tau)&=&\frac{1}{n}\sum_{j=1}^n
\nabla Q_\tau(\X_j)\bigl\{\nabla Q_\tau(
\X_j)\bigr\}^\top
\nonumber\\[-8pt]\\[-8pt]
&&{}+h_n^{-(p+1)}\bigl[\tilde M_n(\tau)\tilde
\Gamma^\top +\tilde\Gamma\tilde M_n^\top(\tau)
\bigr] +o\bigl(n^{-1/2}\bigr),\nonumber
\end{eqnarray}
where $\tilde\Gamma$ is the $p\times s(A)$ matrix consisting of the
second up to the ($p+1)$th rows of
$\Gamma^{-1}$ and
\[
\tilde M_n(\tau)=\frac{1}{n^2}\sum
_{i,j}\frac{\nabla
Q_\tau(\X_j) \X^\top_{ij}(h_{n},A)}{
g(\X_j|\tau)f_{\X}(\X_j)} \bigl[I\bigl\{Y_{i}\le
Q_\tau(\X_{i})\bigr\}-\tau\bigr]I\bigl(|\X_{ij}|\le
h_n\bigr).
\]
The key to the study of the properties of $\hat\Sigma(\tau)$ is $\{
\tilde M_n(\tau)\dvtx \tau\in(0,1)\}$,
which is a typical example of U-processes [\citet{r27}].

To derive the Hoeffding's decomposition of $\tilde M_n(\tau)$, write
$\Z_i=(Y_i,\X_i)$ and define
%
%eA.7 #&#
\begin{eqnarray}\label{fu}
\xi_n(\Z_i,\Z_j;\tau)
&=&   \biggl\{\frac{\nabla Q_\tau
(\X_j)\X^\top_{ij}(h_{n},A)}{g(\X_j|\tau)f_{\X}(\X_j)} \bigl[I\bigl\{Y_{i}\le Q_\tau(
\X_{i})\bigr\}-\tau\bigr]\nonumber
\\
&&\hspace*{5pt}{}+ \frac{\nabla Q_\tau(\X_i) \X^\top
_{ji}(h_{n},A)}{g(\X_i|\tau)f_{\X}(\X_i)} \bigl[I\bigl\{Y_{j}\le
Q_\tau(\X_{j})\bigr\}-\tau\bigr] \biggr\}I\bigl(|
\X_{ij}|\le h_n\bigr),\hspace*{-7pt}\nonumber
\\
\zeta_n(\Z_i;\tau)&=& E_j\bigl[
\xi_n(\Z_i,\Z_j;\tau)\bigr]
\nonumber\\[-4pt]\\[-12pt]
&=& h_n^{p}\bigl[I\bigl\{Y_{i}\le
Q_\tau(\X_{i})\bigr\}-\tau\bigr]\nonumber
\\
&&{}\times \biggl\{\frac{ \nabla
Q_\tau(\X_i) }{g(\X_i|\tau)} \gamma^\top\nonumber
\\
&&\hspace*{18pt}{} +h_n \biggl[\frac{\nabla^2
Q_\tau(\X_i)}{g(\X_i|\tau)}-\frac{ \nabla
Q_\tau(\X_i) \nabla^\top g(\X_i|\tau)}{g^2(\X_i|\tau)} \biggr]
\Gamma_1+O\bigl(h_n^2\bigr) \biggr\},\nonumber
\end{eqnarray}
where
\[
{\gamma}=\int_{[-1,1]^{p}} \t(A)\,d\t,\qquad\Gamma_1=
\int\t \t ^\top(A)\,d\t.
\]
Note that $E[\xi_n(\Z_i,\Z_j;\tau)]=E[\zeta_n(\Z_i;\tau)]=0$.
Therefore, we have
\begin{eqnarray*}
\tilde M_n(\tau)&=&\frac{1}{n^2}\sum
_{i<j}\xi_n(\Z_i,\Z _j;
\tau)=U_n(\tau) %\label{wen}
+\frac{1}{n}\sum
_{i}\zeta_n(\Z_i;\tau),
\end{eqnarray*}
where $U_n(\tau)$ is its Hoeffding's decomposition
%
%eA.8 #&#
\begin{equation}
U_n(\tau)=\frac{1}{n^2}\sum_{i<j}
\xi_n(\Z_i,\Z_j;\tau) -\frac{1}{n}\sum
_{i}\zeta_n(\Z_i;\tau).
\label{leverage}
\end{equation}
To decide the tail properties of
$\sup\{|U_n(\tau)|\dvtx \tau\in[\delta^*,1-\delta^*]\}$, first note
that according to
Lemma 2.13 of \citet{r28} [reproduced as (C1) at the end
of this section] and Corollary \ref{garden},
$\{\xi_n(\Z_i,\Z_j;\tau)\dvtx \tau\in[\delta^*,1-\delta^*]\}$ is
Euclidean for a constant envelope, or in
\citet{r1} term, a uniformly bounded V--C subgraph class.
Applying Proposition 4 in \citet{r1} to $U_n(\tau)$, we conclude that
there exists some finite $c_2>0$, such that for any $\epsilon>0$,
\[
P\Bigl\{n^{1/2}\sup_{\tau\in[\delta^*,1-\delta^*]}\bigl|U_n(\tau )\bigr|\ge
h^{p+1}_n\epsilon\Bigr\}\le 2\exp\bigl\{-c_2
\epsilon n^{1/2}h_n^{-1}\bigr\}.
\]
By an application of the Borel--Cantelli lemma, we have
\[
\sup_{\tau\in[\delta^*,1-\delta^*]}\bigl|U_n(\tau )\bigr|=o\bigl(n^{-1/2}h^{p+1}_n
\bigr)\qquad\mbox{a.s.} %\label{ning}
\]
This together with (\ref{ling}), (\ref{fu}), (\ref{leverage}) and
the facts
that $\tilde\Gamma\gamma={\mathbf0}$, $\tilde\Gamma\Gamma_1={\mathbf I}_p$
implies that with probability one,
\begin{eqnarray*}
\hat\Sigma(\tau)&=&\frac{1}{n}\sum_{i=1}^n
\nabla Q_\tau(\X_i)\bigl\{\nabla Q_\tau(
\X_i)\bigr\}^\top
\\
&&{} +\frac{1}{n}\sum
_{i=1}^n \frac{[I\{Y_{i}\le
Q_\tau(\X_{i})\}-\tau]}{g^2(\X_i|\tau)}
\\
&&\hspace*{36pt}{}\times \bigl[2g(\X_i|\tau)\nabla^2 Q_\tau(
\X_i)-\nabla Q_\tau(\X_i)
\nabla^\top g(\X_i|\tau)
\\
&&\hspace*{141pt}{} -\nabla g(\X_i|\tau)\nabla^\top Q_\tau(\X_i) \bigr] +o
\bigl(n^{-1/2}\bigr), %\label{nano}
\end{eqnarray*}
where the term $o(n^{-1/2})$ is uniform in $\tau\in[\delta
^*,1-\delta^*]$.
Consequently, we have
%
%eA.9 #&#
\begin{eqnarray} \label{orwell}
\qquad\hat\Sigma_{\mathrm{T}}&=&\int_{\delta^*}^{1-\delta^*}
\hat\Sigma(\tau)\,d\tau
\nonumber\\[-9pt]\\[-9pt]
&=&\Sigma _{\mathrm{T}}+ \frac{1}{n}\sum
_{i=1}^n \Sigma^{(1)}(\X_i)
+\frac{1}{n}\sum_{i=1}^n
\Sigma^{(2)}(\X _i,Y_i)+o\bigl(n^{-1/2}
\bigr),\qquad\mbox{a.s.},\nonumber
\end{eqnarray}
where $\Sigma^{(1)}(\cdot)$ and $\Sigma^{(2)}(\cdot)$ are two
symmetric random matrices with properties such that
\begin{eqnarray*}
E\bigl[\Sigma^{(1)}(\X)\bigr]&=&{\mathbf0},\qquad E\bigl[
\Sigma^{(2)}(\X,Y)\bigr]={\mathbf0},
\\
\Sigma^{(1)}(\X)\Pi&=&{\mathbf0},\qquad \Sigma^{(2)}(\X,Y)\Pi={\mathbf0},
\end{eqnarray*}
with $\Pi={\mathbf I}-{\mathbf B}_0({\mathbf B}_0^\top{\mathbf B}_0)^{-1}{\mathbf
B}^\top_0$, the projection matrix such that $\Pi{\mathbf B}_0={\mathbf
B}_0^\top\Pi={\mathbf0}$. An\vspace*{1pt} application of Lemma A.1 in \citet{r23} to\vspace*{1pt}
the right-hand side of (\ref{orwell}) with $\Sigma_{\mathrm{T}}$,
$n^{-1/2}$, $\hat\Sigma_{\mathrm{T}}$ and $n^{-1/2}\sum_i
\{\Sigma^{(1)}(\X_i)+\Sigma^{(2)}(\X_i,Y_i)\}$ acting as $T$, $w^2$,
$T(w)$ and $T^{(2)}$ therein, respectively, we have with probability one,
\begin{eqnarray*}
\tilde\lambda_{p-q}(\hat\Sigma_{\mathrm{T}}) &=&\frac{n^{-1/2}}{p-q}
\sum_{i}\operatorname{trace}\bigl(\bigl[\Sigma^{(1)}(
\X_i)+ \Sigma^{(2)}(\X_i,Y_i)
\bigr]\Pi\bigr)+o\bigl(n^{-1/2}\bigr)
\\
&=&o\bigl(n^{-1/2}\bigr).
\end{eqnarray*}

We now move on to derive the asymptotic properties of the first $q$
eigenvalues and eigenvectors of $\hat\Sigma$.
First note that the three classes of functions, namely $\{\nabla
Q_\tau(\X_i)\{\nabla Q_\tau(\X_i)\}^\top,  \tau\in[\delta
^*,1-\delta^*]\}$,
$\{g[(\X_i|\tau)]^{-2}[I\{Y_{i}\le Q_\tau(\X_{i})\}-\tau], \tau
\in[\delta^*,1-\delta^*]\}$,
and $ \{g(\X_i|\tau)\nabla^2
Q_\tau(\X_i)-\nabla
Q_\tau(\X_i) \nabla^\top\hspace*{-1pt} g(\X_i|\tau)-\nabla
g(\X_i|\tau) \nabla^\top\!\* Q_\tau(\X_i),\tau\in[\delta^*,1-\delta
^*]\}$ are,
according to Corollary \ref{garden}, all Euclidean for a constant envelope.
Therefore, the collection of random matrices
$\{\hat\Sigma(\tau)\dvtx \tau\in[\delta^*,1-\delta^*]\}$ are \emph{Glivenko--Cantelli} as
well as \emph{Donsker} [\citet{r32}].

By \emph{Glivenko--Cantelli}, we mean that
\[
\sup_{\tau\in[\delta^*,1-\delta^*]}\bigl|\operatorname{Vech}\bigl(\hat \Sigma(\tau)\bigr)-
\operatorname{Vech}\bigl(\Sigma(\tau)\bigr)\bigr|\to0\qquad\mbox{a.s.},
\]
from which we can conclude that
\[
\operatorname{Vech}(\hat\Sigma_{\mathrm{T}})-\operatorname{Vech}(
\Sigma_{\mathrm{T}})\to0\qquad\mbox{a.s.}%\label{cantelli}
\]
which in turn implies that [Lemma 3.1, \citet{r3}],
\[
\beta_k(\hat\Sigma_{\mathrm{T}})-\beta_k(
\Sigma_{\mathrm{T}})\to0\qquad (k=1,\ldots,q)\mbox{ a.s.}
\]

By \emph{Donsker}, we mean that
\[
\sqrt{n}\bigl\{\operatorname{Vech}\bigl(\hat\Sigma(\tau)\bigr)-
\operatorname{Vech}\bigl(\Sigma(\tau )\bigr)\bigr\} \stackrel{d} {\to} \mathbb{G}
\qquad\mbox{in }\ell^{\infty
}\bigl(\bigl[\delta^*,1-\delta^*\bigr]\bigr),
\]
where $\ell^{\infty}([\delta^*,1-\delta^*])$ stands for the space
of all uniformly bounded multivariate real functions
from $[\delta^*,1-\delta^*]$ to $R^{p(p+1)/2}$ equipped with the
supremum norm, and the limit $\mathbb{G}$
is a zero-mean $p(p+1)/2$-dimensional Gaussian process on $[\delta
^*,1-\delta^*]$, such that
for any given $\tau_1,\tau_2\in[\delta^*,1-\delta^*]$, the
covariance matrix $E[\mathbb{G}(\tau_1)\mathbb{G}(\tau_2)]$
has its $(k,l)$th element given by the covariance between
\begin{eqnarray*}
&& \nabla Q^{[v(k,1)]}_{\tau_1}(\X)\nabla Q^{[v(k,2)]}_{\tau_1}(\X)
\\
&&\qquad{} + \frac{[I\{Y_{i}\le
Q_{\tau_1}(\X_{i})\}-\tau_1]}{g^2(\X_i|\tau_1)}
\bigl[2g(\X_i|\tau _1)
\nabla^2_{[v(k,1),v(k,2)]} Q_{\tau_1}(\X_i)
\\
&&\hspace*{143pt}{}-
\nabla Q^{[v(k,1)]}_{\tau_1}(\X_i)
\nabla^{[v(k,2)]} g(\X_i|\tau _1)
\\
&&\hspace*{143pt}{} -
\nabla^{[v(k,1)]} g(\X_i|\tau_1)
\nabla^{[v(k,2)]} Q_{\tau_1}(\X_i) \bigr]
\end{eqnarray*}
and
\begin{eqnarray*}
&& \nabla Q^{[[v(l,1)]]}_{\tau_2}(\X)\nabla Q^{[v(l,2)]}_{\tau_2}(\X)
\\
&&\qquad{} + \frac{[I\{Y_{i}\le
Q_{\tau_2}(\X_{i})\}-\tau_2]}{g^2(\X_i|\tau_2)} \bigl[2g(\X_i|\tau _2)
\nabla^2_{[v(l,1),v(k,2)]} Q_{\tau_2}(\X_i)
\\
&&\hspace*{143pt}{}-\nabla Q^{[v(l,1)]}_{\tau_2}(\X_i)
\nabla^{[v(l,2)]} g(\X_i|\tau _2)
\\
&&\hspace*{143pt}{}-
\nabla^{[v(l,1)]} g(\X_i|\tau_2)
\nabla^{[v(l,2)]} Q_{\tau_2}(\X_i) \bigr];
\end{eqnarray*}
equation (\ref{donsker}) thus follows by appealing to the
continuous-mapping theorem.

The proof of (\ref{evalue}) and (\ref{evector}), that is, the
asymptotic normality of
the eigenvalues and eigenvectors of $\hat\Sigma$,
can be done in exactly the same manner as in Theorem 2.2 of \citet{r40},
which by an application of the perturbation theory [\citet{r30}, \citet{r14}],
relates the asymptotic normality of a random matrix to
that of its eigenvalues and eigenvectors.
\end{pf*}

To prepare for the proof of Lemma \ref{T1}, we need to introduce
more notation and some related results.
For any given $\mathbf{x}\in\mathcal{D}$,
let $DX_{n}(\mathbf{x})$ be the $N_n(\mathbf{x})\times s(A)$ matrix
with rows
given by the transposition of $\X
_{i\mathbf{x}}(h_{n},A)$, $i\in
S_{n}(\mathbf{x})$, and $VY_{n}(\mathbf{x})$ be the $N_n(\mathbf
{x})\times1$
vector whose components are $Y_{i}$, $i\in S_{n}(\mathbf{x})$.

For any subset $\mathbf{h}\subset S_{n}(\mathbf{x})$, denote by
$DX_{n}(\mathbf{x},\mathbf{h})$ and $VY_{n}(\mathbf{x},\mathbf
{h})$, the sub-matrix (vector) of $DX_{n}(\mathbf{x})$ and
$VY_{n}(\mathbf{x})$, respectively,
with indices of rows given by $\mathbf{h}$. Further define
\[
\H_n(\mathbf{x})=\bigl\{\mathbf{h}\dvtx \mathbf{h}\subset
S_{n}(\mathbf{x} ), \sharp(\mathbf{h})=s(A), DX_{n}(
\mathbf{x},\mathbf{h}) \mbox{ is of full rank}\bigr\}.
\]
Suppose $DX_{n}(\mathbf{x})$ of rank $=s(A)$, $\H_n(\mathbf{x})$ is
thus nonempty.
The following two facts concern the uniqueness of $\hat\c_n(\mathbf
{x};\tau)$ and its ``matrix form'' of, for
any given $\mathbf{x}\in\mathcal{D}$ and $\tau\in(0,1)$.
They are essentially restatements of Theorems 3.1 and 3.2 in \citet{r16};
see also Facts 6.3 and 6.4 in \citet{r4}.
\begin{longlist}[(B3)]
\item[(B1)] There exist positive constants $c_1$ and $c_2$, such that
\[
P(A_n)=1\qquad\mbox{where }A_n=\bigl
\{c_1nh_n^d\le N_n(\mathbf{x})\le
c_2nh_n^d\mbox{ for all }\mathbf{x}\in
\mathcal{D}\bigr\}.
\]
This follows easily from (\ref{ellen}).
\item[(B2)] There exists a
$\mathbf{h}\in\H_n(\mathbf{x})$, such that (\ref{a8}) has at least one
solution of the form
\[
\hat{\c}_{n}(\mathbf{x};\tau)=\bigl[DX_{n}(\mathbf{x},
\mathbf{h})\bigr]^{-1} VY_{n}(\mathbf{x},\mathbf{h}).
\]

\item[(B3)] For $\mathbf{h}\in\H_n(\mathbf{x})$,
let $\hat\c_n(\mathbf{x};\tau)=[DX_{n}(\mathbf{x},\mathbf{h})]^{-1}
VY_{n}(\mathbf{x},\mathbf{h})$ and define
\[
L_{n}(\mathbf{h};\mathbf{x},\tau) =\bigl[DX_{n}(
\mathbf{x},\mathbf{h} )\bigr]^{-1}\sum_{i\in\bar{\mathbf{h}}}
\bigl[I\bigl\{ Y_{i}< \X _{i\mathbf{x}}^\top(h_{n},A)
\hat{\c}_{n}(\mathbf{x};\tau) \bigr\}-\tau \bigr]\X
_{i,\mathbf{x}}(h_{n},A),
\]
where $\bar{\mathbf{h}}$ is the relative complement
of $\mathbf{h}$ with respect to $S_{n}(\mathbf{x})$. Then $\hat{\c
}_{n}(\mathbf{x};\tau)$ is a
unique solution to (\ref{a8}) if and only if $L_{n}(\mathbf
{h};\mathbf{x},\tau)\in
(\tau-1,\tau)^{ s(A)}$. Further, if $\hat{\c}_{n}(\mathbf{x};\tau
)$ is a
solution (not necessarily unique) to (\ref{a8}), we must have
$L_{n}(\mathbf{h};\mathbf{x},\tau)\in
[\tau-1,\tau]^{ s(A)}$.
\end{longlist}
To facilitate the use of the conditioning arguments at various places
in the proofs,
for any $\X_j$, $j=1,\ldots,n$, we exclude $\X_j$ from the previously
defined $S_n(\X_j)$;
instead we define $S_n(\X_j)=\{i\dvtx 1\le i\le n, i\ne
j,  |\X_{ij}|\le h_n\}$ and $N_n(\X_j)=\sharp(S_n(\X_j))$.

The proof of Lemma \ref{T1} will be built upon the following slightly
weaker result.

%
%leA.1 #&#
\begin{Lemma}\label{lemma2}Let $\delta_{n}=(nh_n^p/\log n)^{-1/2}$.
Suppose conditions in Lemma \ref{T1} hold.
Then
\[
\sup_{1\le j\le n,\tau\in[\delta^*,1-\delta^*]}\bigl|\hat\c _{n}(\X_j;\tau)-
\c_{n}(\X_j;\tau)\bigr|=O(\delta_n)\qquad
\mbox{a.s.}
\]
\end{Lemma}

\begin{pf} %{Proof of Lemma \ref{lemma2}}
For any given positive
constant $K_{1}$ and a generic $\mathbf{x}\in\mathcal{D}$, let
$U_{n}$ be the event defined as
%
%eA.10 #&#
\begin{equation}
U_{n}= \Bigl\{\sup_{\tau\in[\delta^*,1-\delta^*]} \bigl|\hat\c_{n}(
\mathbf{x};\tau)-\c_{n}(\mathbf{x};\tau)\bigr|\geq K_1
\delta_n \Bigr\}.\label{else}
\end{equation}
In view of the fact that $P(A_n)=1$, the assertion in Lemma \ref
{lemma2} will
follow from an application of the Borel--Cantelli lemma,
if we can show that there exists some $K_{1}>0$, such that
%
%eA.11 #&#
\begin{equation}
\sum_{n}nP(U_{n}\cap A_n)<
\infty.\label{b1}
\end{equation}
We now try to get an upper bound for $P(U_{n}\cap A_n)$. To this end, for
given $\tau\in[\delta^*,1-\delta^*], \mathbf{x}\in\mathcal{D}$
and $\c\in R^{s(A)}$, define
\begin{eqnarray}
\nonumber
Z_{ni}(\c|\mathbf{x},\tau)&=& \bigl[I\bigl\{Y_i<
\c^\top \X_{i\mathbf{x}}(h_{n},A) \bigr\}-\tau \bigr]\X
_{i,\mathbf{x}}(h_{n},A). % \label{julia}
\end{eqnarray}
Based on (B2) and (B3), there exists some positive constant $K_2$,
which depends only on $s(A)$ such that $U_{n}\cap A_n$ is contained in
the event
%
%eA.12 #&#
\begin{eqnarray}\label{b2}\quad
&& \biggl\{\mbox{there exists some }\tau\in\bigl[\delta^*,1-\delta^*
\bigr] \mbox{ and }\mathbf{h}\in\H_n(\mathbf{x}),\mbox{ such that for}\hspace*{-25pt}\nonumber
\\
&&\hspace*{6pt}\hat{\c}_{n}(\mathbf{x};\tau)=\bigl[DX_{n}(
\mathbf{x},\mathbf{h})\bigr]^{-1}VY_{n}(\mathbf{x},
\mathbf{h}),\mbox{ we have }\hspace*{-25pt}
\\
&&\hspace*{9pt} \Biggl|\sum_{i\in\bar{\mathbf{h}}}Z_{ni}
\bigl(\hat{\c}_{n}(\mathbf {x};\tau)|\mathbf{x},\tau\bigr)\Biggr|\le
K_{2}\mbox{ and }\bigl|\hat{\c}_{n}(\mathbf{x};\tau)-{\c}_{n}(\mathbf
{x};\tau)\bigr|\geq K_{1}\delta_n \biggr\}\cap A_n.\hspace*{-25pt}\nonumber
\end{eqnarray}
Choose large enough $K_1$
such that we can apply Proposition \ref{a10} to conclude that
there exist some $\epsilon_{1}>0$, and $K_3>0$,
such that, for all $\tau\in[\delta^*,1-\delta^*]$,
\[
E\bigl[Z_{ni}\bigl(\hat{\c}_{n}(\mathbf{x};\tau)|
\mathbf{x},\tau\bigr)\bigr] \ge \min\{\epsilon_{1},K_3K_1
\delta_n\},
\]
and consequently as a result of $A_n$ and the fact that $\sharp(\bar
{\mathbf{h}})=N_n(\mathbf{x})-s(A)$, we have
%
%eA.13 #&#
\begin{eqnarray}\label{a9}
&&\biggl\{ \biggl|\sum_{i\in\bar{\mathbf{h}}}Z_{ni}
\bigl(\hat{\c}_{n}(\mathbf {x};\tau)|\mathbf{x},\tau\bigr) \biggr|\le
K_{2} \biggr\}
\nonumber\hspace*{-30pt}\\[-4pt]\\[-12pt]
&&\qquad \subseteq \biggl\{ \biggl|\sum_{i\in\bar{\mathbf{h}}}\bigl
\{Z_{ni}\bigl(\hat{\c}_{n}(\mathbf {x};\tau)|\mathbf{x},\tau
\bigr) -E\bigl[Z_{ni}\bigl(\hat{\c}_{n}(\mathbf{x};\tau)|
\mathbf{x},\tau\bigr)\bigr]\bigr\} \biggr|\ge c_1^*K_1nh_n^p
\delta_n \biggr\}\nonumber\hspace*{-30pt}
\end{eqnarray}
for some $c_1^*>0$.

Next, note that given the set $S_{n}(\mathbf{x})$, $\mathbf{h}\subset
S_{n}(\mathbf{x})$, and $(\X_i,Y_i)$ for $i\in\mathbf{h}$, and thus
$\hat{\c}_{n}(\mathbf{x};\tau)=[DX_{n}(\mathbf{x},\mathbf{h})]^{-1}VY_{n}(\mathbf{x},\mathbf{h})$ is also fixed,
the random vectors $\{Z_{ni}(\hat{\c}_{n}(\mathbf{x};\tau)|\mathbf
{x},\tau), i\in\bar{\mathbf{h}}\}$ are conditionally\vspace*{1pt} i.i.d.
This together with (\ref{b2}), (\ref{a9}) and the fact that $\sharp
(\H_n(\mathbf{x}))$ is of order $(nh_n^p)^{s(A)}$,
implies there exists some $c_2^*>0$, such that
%
%eA.14 #&#
\begin{eqnarray}\label{derailed}
\quad&&  P(U_{n}\cap A_n)\nonumber
\\
&&\qquad \le c_2^*\bigl(nh_n^p\bigr)^{s(A)}
\\
&&\qquad\quad{}\times P \biggl\{\mathop{\sup_{\tau\in[\delta^*,1-\delta
^*],}}_{\c\in R^{s(A)}} \biggl|
\sum_{i\in\bar{\mathbf{h}}}\bigl\{Z_{ni}(\c|\mathbf{x},\tau
)-E\bigl[Z_{ni}(\c|\mathbf{x},\tau)\bigr]\bigr\} \biggr|\ge
c_1^*K_1nh_n^p
\delta_n \biggr\}.\hspace*{-15pt}\nonumber
\end{eqnarray}
To find a bound for the probability on the right-hand side above, first
note that
according to Lemma 22(ii) in \citet{r27},
$\{Z_{ni}(\c|\mathbf{x},\tau)\dvtx \tau\in[\delta^*,1-\delta^*],\c
\in R^{s(A)}\}$ is contained in
a Euclidean class for a constant envelope, since
$Y_i-\c^\top\X_{i\mathbf{x}}(h_{n},A)
=[\X^\top_{i\mathbf{x}}(h_{n},A),Y_i]*(\c^\top,-1)^\top$ and the indicator
function $I(\cdot< 0)$ is of bounded variation. As
$E|Z_{ni}(\c|\mathbf{x},\tau)Z^\tau_{ni}(\c|\mathbf{x},\tau
)|^2=O(1)$ uniformly in
$\tau\in[\delta^*,1-\delta^*],\c\in R^{s(A)}$, through similar
arguments used in the proof of
Theorem~2.37 in \citeauthor{r29} [(\citeyear{r29}), page~34], we have that
\[
P \biggl\{\mathop{\sup_{\tau\in[\delta^*,1-\delta^*],}}_{\c\in
R^{s(A)}} \biggl|\sum
_{i\in\bar{\mathbf{h}}}\bigl\{Z_{ni}(\c|\mathbf{x},\tau )-E
\bigl[Z_{ni}(\c|\mathbf{x},\tau)\bigr]\bigr\} \biggr|\ge c_1^*K_1nh_n^p
\delta_n \biggr\}=o\bigl(n^{-a}\bigr),
\]
for any $a>0$. This together with (\ref{derailed}) leads to (\ref{b1}).
\end{pf}

For any $\mathbf{x}\in\mathcal{D}$, let $\omega_{h_n}(\t|\mathbf
{x})$ be the conditional probability density function
of $(\X_i-\mathbf{x})/h_n$ given $i\in S_n(\mathbf{x})$. Note that
it converges to the uniform density on $[-1,1]^p$
uniformly in $\t\in[-1,1]^p$ and $\mathbf{x}\in\mathcal{D}$.

\begin{pf*}{Proof of Lemma \ref{T1}}
For any given $\tau\in[\delta^*,1-\delta^*]$,
$\mathbf{x}\in\mathcal{D}$, and $\X\in S_n(\mathbf{x})$, write
\[
\hat Q_{n}(\X,\mathbf{x};\tau)=\bigl[(\X-\mathbf{x})
(h_n,A)\bigr]^\top\hat \c_{n}(\mathbf{x};\tau).
\]
The proof consists of the following steps.
\begin{longlist}[\textit{Step} 2:]
\item[\textit{Step} 1:] For any given $\tau\in[0,1]$,
$\c\in R^{s(A)}$ and
$\mathbf{x}\in R^{p}$, define
\begin{eqnarray*}
\tilde{H}_{n}(\c;\mathbf{x})&=&E\bigl[I\bigl\{Y_{i}<
c^\top\X _{i\mathbf{x}}(h_{n},A)\bigr\}\X
_{i\mathbf{x}}(h_{n},A)|i\in S_n(\mathbf{x})\bigr]
\\
&=&\int_{[-1,1]^{p}}F\bigl(\c^\top\t(A)|
\mathbf{x}+h_n\t\bigr)\t(A)\omega _{h_n}(\t|\mathbf{x})\,d
\t,
\\
R_{n}^{(1)}(\tilde\c,\c|\mathbf{x},\tau)&=&
\tilde{H}_{n}(\mathbf{x},\tilde\c)- \tilde{H}_{n}(
\mathbf{x},\c) -\Sigma_n(\mathbf{x};\tau) (\tilde\c-\c).
\end{eqnarray*}
Therefore, under assumptions (A2) and (A3),
%
%eA.15 #&#
%eA.16 #&#
\begin{eqnarray}
\nonumber
&&R_{n}^{(1)}\bigl(\hat\c_{n}(
\mathbf{x};\tau),\c _{n}(\mathbf{x};\tau)|\mathbf{x},\tau\bigr)\hspace*{-20pt}
\\
&&\qquad=\tilde{H}_{n}\bigl(\mathbf{x},\hat\c_{n}(
\mathbf{x};\tau)\bigr)- \tilde{H}_{n}\bigl(\mathbf{x},
\c_{n}(\mathbf{x};\tau)\bigr) -\Sigma_n(\mathbf{x};\tau)
\bigl[\hat\c_{n}(\mathbf{x};\tau)-\c _{n}(\mathbf{x};\tau)
\bigr]\hspace*{-20pt}\label{monica}
\\
\nonumber
&&\qquad=\int_{[-1,1]^{p}}\bigl[F\bigl(\hat
Q_{n}(\mathbf {x}+h_n\t,\mathbf{x};\tau)|
\mathbf{x}+h_n\t\bigr)\hspace*{-20pt}
\\
&&\hspace*{35pt}\quad\qquad{} -F\bigl(Q_{n}(\mathbf{x}+h_n
\t,\mathbf{x};\tau)|\mathbf{x}+h_n\t\bigr)\hspace*{-20pt}\nonumber
\\
\nonumber
&&\hspace*{35pt}\quad\qquad{}-g(\mathbf{x}+h_n\t|\tau)\t(A)
\t^\top(A)\bigl\{ \hat\c_{n}(\mathbf{x};\tau)-
\c_{n}(\mathbf{x};\tau)\bigr\} \bigr]w_{h_n}(\t|\mathbf{x})
\,d\t\hspace*{-20pt}
\\
&&\qquad=O\bigl(\delta_n^{1+s_3}\bigr)=O\bigl\{
\bigl[n^{(1-\kappa p)}/\log n\bigr]^{-3/4}\bigr\}\qquad (\mbox {if
}s_{3}\ge1/2),\hspace*{-20pt} \label{phoebe}
\end{eqnarray}
uniformly in $\tau\in[\delta^*,1-\delta^*]$, where (\ref{phoebe})
follows from Lemma \ref{lemma2} and the facts that
$\hat Q_{n}(\mathbf{x}+h_n\t,\mathbf{x};\tau)-Q_{n}(\mathbf
{x}+h_n\t,\mathbf{x};\tau)=\{\t(A)\}^\top[\hat\c_{n}(\mathbf
{x};\tau)-\c_{n}(\mathbf{x};\tau)]$
and $Q_{n}(\mathbf{x}+h_n\t,\mathbf{x};\tau)-Q_\tau(\mathbf
{x}+h_n\t)=O(h_n^{s_2})=o(\delta_n)$.

\item[\textit{Step} 2:] For any given $\tau\in(0,1)$,
$\mathbf{x}\in R^{p}$ and $\mathbf{h}\in H_{n}(\mathbf{x})$, define
\begin{eqnarray*}
\chi_{n}(\mathbf{x};\tau) & =&\sum_{i\in S_{n}(\mathbf
{x})}
\bigl[\X_{i\mathbf{x}}(h_{n},A)I\bigl\{Y_{i}%
\le\hat Q_{n}(\X_{i},\mathbf{x};\tau)\bigr\}- \tilde
H_{n}\bigl(\hat\c _{n}(\mathbf{x};\tau);\mathbf{x}\bigr)
\bigr]
\\
&&{} -\sum_{i\in S_{n}(\mathbf{x})} \bigl[\X_{i\mathbf
{x}}(h_{n},A)I
\bigl\{Y_{i}\le Q_{n}(\X%
_{i},
\mathbf{x};\tau)\bigr\}- \tilde H_{n}\bigl( \c _{n}(
\mathbf{x};\tau),\mathbf{x}\bigr)\bigr],
\\
\hat{\c}^{\mathbf{h}}_{n}(\mathbf{x};\tau) & =&
\bigl[DX_{n}(\mathbf{x},\mathbf{h})\bigr]^{-1}VY_{n}(
\mathbf{x},\mathbf{h}),
\\
\hat Q_{n}^{\mathbf{h}}(\X_{i},\mathbf{x};\tau)&=&\bigl\{\hat{\c}^{\mathbf{h}}
_{n}(\mathbf{x};\tau)\bigr\}^\top\X_{i\mathbf{x}}(h_{n},A),
\end{eqnarray*}
and for any $\c_1$, $\c_2\in R^{s(A)}$, define
\begin{eqnarray*}
\chi_{n}^{\mathbf{h}}(\c_1,\c_2;\mathbf{x})
& =&\sum_{i\in\bar{\mathbf{h}}} \bigl[ \X_{i\mathbf{x}}(h
_{n},A)I\bigl\{Y_{i}\le\c_1^\top
\X_{i\mathbf{x}}(h_{n},A)\bigr\}- \tilde H_{n}(
\c_1;\mathbf{x})\bigr]
\\
&&{} -\sum_{i\in\bar{\mathbf{h}}} \bigl[\X_{i\mathbf
{x}}(h_{n},A)I
\bigl\{Y_{i}\le\c_2^\top\X_{i\mathbf{x}}(h_{n},A)
\bigr\}- \tilde H_{n}(\c_2;\mathbf{x})\bigr].
\end{eqnarray*}
For any given $K_{3}>0$, consider the corresponding event
\begin{eqnarray*}
W_{n}(\mathbf{x})&=& \Bigl\{\sup_{\tau\in[\delta^*,1-\delta
^*]}\bigl|
\chi_{n}(\mathbf{x};\tau)\bigr|\ge K_{3}[\log
n]^{3/4}%
n^{(1-\kappa p)/4} \Bigr\}.
\end{eqnarray*}
Then in view of definition of the events $A_{n}$,
$U_{n}(\mathbf{x})$ of (\ref{else}) and (B2), the event
$W_{n}(\mathbf{x})\cap A_{n}%
\cap\overline{U_{n}{(\mathbf{x})}}$ [$\overline{U_{n}{(\mathbf{x})}}$
is the complement of $U_{n}(\mathbf{x})$] is contained in the event
\begin{eqnarray*}
&& \bigl\{\mbox{for some }\tau\in\bigl[\delta^*,1-\delta^*\bigr]\mbox{ and }
\mathbf{h}\in\H_{n}(\mathbf{x}),
\\
&&\hspace*{4pt}\bigl|\chi_{n}^{\mathbf{h}%
}\bigl(\hat{\c}^{\mathbf{h}}_{n}(\mathbf{x};\tau),\c_n(\mathbf {x};\tau); \mathbf{x}\bigr)\bigr|
\ge K_{4}[\log n]^{3/4}n^{(1-\kappa p)/4} \mbox{ and}
\\
&&\hspace*{132pt}\bigl|\hat{\c}^{\mathbf{h}}_{n}(\mathbf{x})-\c%
_{n}(\mathbf{x};\tau)\bigr|\le K_{1}\delta_n \bigr\}\cap A_{n}
\end{eqnarray*}
for large enough $n$, where $K_{4}=K_{3}/2$ and for which we have
implicitly used the
facts that $\sharp
(\mathbf{h})=p$ and $[\log n]^{3/4}n^{(1-\kappa p)/4}\to\infty$ as
$n\to\infty$. Again, since $\sharp( H_{n}(\mathbf{x}))$ is of order
$n^{(1-\kappa p)n(A)}$ uniformly in
$\mathbf{x}\in\mathcal{D}$, there exists some constant $c_3>0$,
such that $P(W_{n}(\mathbf{x})\cap A_{n}\cap\overline{U_{n}(\mathbf
{x})})$ is
bounded by $c_3n^{(1-\kappa p)n(A)}$ multiplied by the probability of
the following event:
%
%eA.17 #&#
\begin{equation}
\Bigl\{\mathop{\sup_{\c_1,\c_2\in R^{s(A)};}}_{|\c_1-\c_2|\le K_{1}
\delta_n}\bigl|
\chi_{n}^{\mathbf{h}}(\c_1,\c_2; \mathbf{x})\bigr|
\ge K_{4}[\log n]^{3/4}n^{(1-\kappa p)/4} \Bigr\}\cap
A_{n}. \label{b3}
\end{equation}
To find a bound for the probability of even (\ref{b3}), first note
that according to
Lemma 22(ii) in \citet{r27} and Lemma 2.14(i) in \citet{r28},
the class of all functions on $R^{s(A)+1}$
of the form
\[
\bigl(Y_i,\X_{i\mathbf{x}}(h _{n},A)\bigr)\to
\X_{i\mathbf{x}}(h _{n},A)\bigl[I\bigl\{Y_{i}\le
\c_1^\top\X_{i\mathbf{x}}(h_{n},A)\bigr\}-I
\bigl\{ Y_{i}\le\c_2^\top\X_{i\mathbf{x}}(h_{n},A)
\bigr\}\bigr]
\]
%
% with
$\c_1,\c_2$ ranging over $ R^{s(A)}$ is again a
Euclidean class for a constant envelope.
Secondly, conditioning on $S_{n}%
(\mathbf{x})$, $\mathbf{h}\in H_{n}(\mathbf{x})$, and observations
$\{(\X_{i},Y_{i})\dvtx\break  i\in\mathbf{h}\} $, the terms in the sum defining $\chi_{n}%
^{\mathbf{h}}(\c_1,\c_2; \mathbf{x})$ are i.i.d. with mean zero, and
variance--covariance matrix with Euclidean norm of order
$O(|\c_1-\c_2|)$. Following the steps in the proof of
Theorem 2.37 in \citeauthor{r29} [(\citeyear{r29}), page~34], we can conclude that
there exist constant $c_{4}>0,c_{5}>0$, such that the probability of
(\ref{b3}) is bounded by
\[
K_4^{c_4}(\log n)^{c_4/2}\exp\bigl(-c_5K_4^2
\log n\bigr)=o\bigl(n^{-\alpha}\bigr)\qquad\mbox{for any }\alpha>0,
\]
if $K_{4}$, or equivalently $K_3$, is chosen to be sufficiently large.
Equivalently, we have
there exists some $K_3$, such that
\begin{eqnarray*}
&&P \Bigl\{\sup_{\tau\in[\delta^*,1-\delta^*]}\bigl|\chi _{n}(\mathbf{x};\tau)\bigr|\ge
K_{3}[\log n]^{3/4}%
n^{(1-\kappa p)/4} \Bigr\}=o
\bigl(n^{-2}\bigr).
\end{eqnarray*}
An application of the Borel--Cantelli lemma leads to
%
%eA.18 #&#
\begin{eqnarray}
&&\sup_{\tau\in[\delta^*,1-\delta^*],j=1,\ldots,n}\bigl|\chi _{n}(\X_j;\tau)\bigr|=O
\bigl\{(\log n)^{3/4}%
n^{(1-\kappa p)/4}\bigr\}\qquad
\mbox{a.s.}\label{kathy}
\end{eqnarray}
\item[\textit{Step} 3:] Combining (\ref{monica}), (\ref{phoebe}) and
(\ref{kathy}), we have
with probability one,
%
%eA.19 #&#
\begin{eqnarray}\label{libya}
\nonumber
&& \frac{1}{N_{n}(\mathbf{x})}\sum_{i\in
S_{n}(\mathbf{x})}
\X_{i\mathbf{x}}(h_{n},A)\bigl[I\bigl\{Y_{i}\le
Q_{n}(\X_{i},\mathbf{x};\tau)\bigr\}-\tau\bigr]
\\
\nonumber
&&\qquad  =-\frac{1}{N_{n}(\mathbf{x})}\chi_{n}^{\mathbf
{h}}(\mathbf{x})-
\tilde H_{n}\bigl(\hat\c_{n}(\mathbf{x};\tau);\mathbf{x}
\bigr) +\tilde H_{n}\bigl(\c_{n}(\mathbf{x};\tau);\mathbf{x}
\bigr)
\\
&&\quad\qquad{} +\frac{1}{N_{n}(\mathbf{x})}\sum_{i\in
S_{n}(\mathbf{x})}
\X_{i\mathbf{x}}(h_{n},A)\bigl[I\bigl\{Y_{i}%
\le\hat Q_{n}(\X_{i},\mathbf{x};\tau)\bigr\}-\tau\bigr]
\\
\nonumber
&&\qquad  = -\Sigma_n(\mathbf{x};\tau)\bigl[\hat
\c_{n}(\mathbf {x};\tau)-\c_{n}(\mathbf{x};\tau)\bigr] +O
\bigl\{\bigl[n^{(1-\kappa p)}/\log n\bigr]^{-3/4}%
\bigr\}
\\
&&\quad\qquad{}+\frac{1}{N_{n}(\mathbf{x})}\sum_{i\in S_{n}(\mathbf{x})}
\X_{ij}(\delta _{n},A)\bigl[I\bigl\{Y_{i}\le\hat
Q_{n}(\X_{i},\mathbf{x};\tau)\bigr\}-\tau
\bigr]\nonumber
\end{eqnarray}
\end{longlist}
uniformly in $\tau\in[\delta^*,1-\delta^*]$ and $\mathbf{x}=\X
_{j}$, $j=1,\ldots,n$.
Note that according
to~(B3), the last term in (\ref{libya}) is of order
$O(n^{\kappa p-1})=o\{[n^{(1-\kappa p)}/\log n]^{-3/4}\}$.
\end{pf*}

%
%prA.2 #&#
\begin{Proposition} \label{a10}There exist some $K_2>0, K_3>0, K_4>0$
such that
for all $\tau\in[\delta^*,1-\delta^*]$,
\[
\biggl|\int_{[-1,1]^p}\bigl\{F\bigl(\c^\top t(A) |
\mathbf{x}+h_n\t\bigr)-\tau\bigr\} t(A)\omega_{h_n}(\t|
\mathbf{x})\,d\t \biggr|\ge \min\bigl\{K_2,K_3\bigl|\c-
\c_{n}(\mathbf{x};\tau)\bigr|\bigr\},
\]
whenever $|\c-\c_{n}(\mathbf{x};\tau)|\ge K_4h_n^{s_2}$.
\end{Proposition}

\begin{pf}%{Proof of Proposition \ref{a10}}
First note that as $\omega
_{h_n}(\t|\mathbf{x})$ converges to
the uniform density on $[-1,1]^p$ uniformly in $\t\in[-1,1]^p$,
$\mathbf{x}\in\mathcal{D}$, we have
%
%eA.20 #&#
\begin{eqnarray}
&&\int_{[-1,1]^p}\bigl\{F\bigl(\c^\top t(A) |
\mathbf{x}+h_n\t\bigr)-\tau\bigr\} t(A)\omega_{h_n}(\t|
\mathbf{x})\,d\t=H_n(\c|\mathbf{x},\tau ) \bigl(1+o(1)\bigr)\nonumber
\\
\eqntext{\displaystyle\mbox{where }H_n(\c|\mathbf{x},\tau)=\int_{[-1,1]^p}
\bigl\{F\bigl(\c^\top t(A) |\mathbf{x}+h_n\t\bigr)-\tau\bigr
\}t(A)\,d\t.}
\end{eqnarray}
The proof is split into the following steps.

\begin{longlist}[\textit{Step} 3:]
\item[\textit{Step} 1:] We show that there exist $M_1>0$ and $\epsilon_1>0$,
such that for all $\tau\in[\delta^*,1-\delta^*]$, and $\mathbf
{x}\in\mathcal{D}$,
$|H_n(\c|\mathbf{x},\tau)|\ge\epsilon_1$, whenever
$|\c-\c_n(\mathbf{x};\tau)|\ge M_1$.

If this is false, there must exist three sequences $\{\tau_{n^*}\}$ in
$ [\delta^*,1-\delta^*]$,
$\{\mathbf{x}_{n^*}\}$ in $\mathcal{D}$ and $\{\c_{n^*}\}$ in $R^{s(A)}$,
such that as $n^*\to\infty$, $|\c_{n^*}-\c_n(\mathbf{x}_{n^*};\tau
_{n^*})|\to\infty$,
but $|H_n(\c_{n^*}|\mathbf{x}_{n^*},\tau_{n^*})|\to0$.
Without loss of generality, suppose there exist some $\tau^*\in
[\delta^*,1-\delta^*]$
and $\mathbf{x}^*\in\mathcal{D}$, such that
as $n^*\to\infty$, $\tau_{n^*}\to\tau^*$, and
$\mathbf{x}_{n^*}\to\mathbf{x}^*$. Further construct the sequence
$\{\Delta_{n^*}\}$ with $\Delta_{n^*}=\c_{n^*}-\c_n(\mathbf
{x}_{n^*};\tau_{n^*})$, and
for which we have, as $n^*\to\infty$,
$|\Delta_{n^*}|\to\infty$, and $\Delta_{n^*}/|\Delta_{n^*}|\to
\Delta^*$,
for some $\Delta^*\in R^{s(A)}$.

Note that for any given $\t\in[-1,1]^p$, $\c_{n^*}^\top\t(A)=
\c_n(\mathbf{x}_{n^*};\tau_{n^*})^\top\t(A)+\Delta_{n^*}^\top\t(A)$,
the first term being finite,
must tend to either $+\infty$
or $-\infty$ depending on whether $\t^\top(A) \Delta^* $ is positive
or negative.
Consequently, due to $F(\cdot|\cdot)$ being continuous in both its arguments,
we have
\begin{eqnarray*}
\lim_{n^*\to\infty}F\bigl(\c_{n^*}^\top\t(A)|
\mathbf {x}_{n^*}+h_n\t\bigr) &=&\lim_{n^*}F
\bigl(\c_{n^*}^\top\t(A)|\mathbf{x}^*+h_n\t\bigr)
\\
&=& F\bigl(+\infty\times\sign\bigl\{\t^\top(A) \Delta^*\bigr\}|
\mathbf{x}^*+h_n\t\bigr),
\end{eqnarray*}
which must tend to either $1$
or $0$ depending on whether $\t^\top(A) \Delta^* $ is positive
or negative, respectively.
As it is trivial to argue that the region $[-1,1]^p\cap\{\t\dvtx \t^\top
(A) \Delta^*=0\}$
must have Lebesgue measure zero,
a simple application of the dominated convergence theorem to $H_n(\c
_{n^*}|\mathbf{x}_{n^*},\tau_{n^*})$ yields
\begin{eqnarray*}
&&\tau^* \int_{[-1,1]^p\cap\{\t\dvtx {\t}^\top(A) \Delta^*<0\}}\t (A)\,d{\t} =\bigl(1-\tau^*\bigr)\int
_{[-1,1]^p\cap\{\t\dvtx {\t}^\top(A) \Delta^*>0\}} \t(A)\,d{\t}.
\end{eqnarray*}
Multiplying either side by $ \Delta^*$, we get
\begin{eqnarray*}
&& \tau^* \int_{[-1,1]^p\cap\{\t\dvtx \t^\top(A) \Delta^*<0\}}\t^\top (A) \Delta^*d{\t}
\\
&&\qquad =
\bigl(1-\tau^*\bigr)\int_{[-1,1]^p\cap\{\t\dvtx \t^\top(A) \Delta^*>0\}} \t^\top(A) \Delta^*d{
\t}.
\end{eqnarray*}
As $0<\tau^*<1$,
the above implies that both regions
$[-1,1]^p\cap\{\t\dvtx \t^\top(A) \Delta^*<0\}$ and
$[-1,1]^p\cap\{\t\dvtx \t^\top(A) \Delta^*>0\}$
must both have Lebesgue measure zero, which cannot be true.

\item[\textit{Step} 2:] For any $\t\in[-1,1]^p$, write $R_n(\t;\tau,\mathbf
{x})=\t^\top(A)\c_n(\mathbf{x};\tau)-Q_\tau(\mathbf{x}+h_n\t)$.
Note that $R_n(\t,\mathbf{x})=O(h_n^{s_2})$ uniformly in $\t\in
[-1,1]^p, \tau\in[\delta^*,1-\delta^*]$
and $\mathbf{x}\in\mathcal{D}\subset R^p$.
For any $\t\in[-1,1]^p$ and $\c\in R^{s(A)}$, define a real valued
function as
\begin{eqnarray*}
&&g_n(\c,{\t}|\mathbf{x},\tau)=\frac{F(\c^\top\t(A)|\mathbf
{x}+h_n\t)-F(\c_n(\mathbf{x};\tau)^\top\t(A)|\mathbf{x}+h_n\t)}{
(\c-\c_n(\mathbf{x};\tau))^\top\t(A) }.
\end{eqnarray*}
In the case where $(\c-\c_n(\mathbf{x};\tau))^\top\t(A) =0$,
$g_n(\c,{\t}|\mathbf{x},\tau)$ can be defined arbitrarily because
the set $\{\t\in[-1,1]^p\dvtx \c^\top\t(A) =0\}$ has Lebesque measure
zero for any nonzero
$\c$.
Write
%
%eA.21 #&#
\begin{eqnarray} \label{a11}
&& H_n(\c|\mathbf{x},\tau)\nonumber
\\
&&\qquad =\int_{[-1,1]^p}
\bigl\{F\bigl(\c^\top\t (A)|\mathbf{x}+h_n\t\bigr)-F\bigl(
\c_n(\mathbf{x};\tau)^\top\t(A)|\mathbf {x}+h_n
\t\bigr)\bigr\}\t(A)\,d\t\nonumber
\\
&&\quad\qquad{}+\int_{[-1,1]^p}\bigl\{F\bigl(\c_n(
\mathbf{x};\tau )^\top\t(A)|\mathbf{x}+h_n\t\bigr)
\nonumber\\[-8pt]\\[-8pt]
&&\hspace*{83pt}{} -F\bigl(Q_\tau(\mathbf{x}+h_n\t)|\mathbf {x}+h_n\t
\bigr)\bigr\}\t(A)\,d\t\nonumber
\\
&&\qquad = \biggl[\int_{[-1,1]^p}g_n(\c,{\t}|
\mathbf{x},\tau)\t(A)\bigl\{\t(A)\bigr\}^\top \,d\mathbf{t} \biggr]\bigl(\c-
\c_n(\mathbf{x};\tau)\bigr)\nonumber
\\
&&\quad\qquad{}+\int_{[-1,1]^p} f_{Y|\X}\bigl(Q_\tau(
\mathbf{x}+h_n\t)+\xi_1R_n(\t;\tau,\mathbf
{x})|\mathbf{x}+h_n\t\bigr) R_n(\t;\tau,\mathbf{x})\t(A)
\,d\t,\nonumber
\end{eqnarray}
where $\xi_1$ lies between $0$ and $1$, depending
on $\t$, $\tau$ and $\mathbf{x}$.

\item[\textit{Step} 3:] By the Cauchy inequality, we have regarding the second term
on the right-hand side of (\ref{a11}),
%
%eA.22 #&#
\begin{eqnarray}\label{a13}
&& \biggl|\int_{[-1,1]^p} \bigl\{f_{Y|\X}
\bigl(Q_\tau(\mathbf{x}+h_n\t)+\xi_1R_n(
\t;\tau,\mathbf {x})|\mathbf{x}+h_n\t\bigr)\bigr\}R_n(\t,
\mathbf{x})\t(A)\,d\t \biggr|^2
\nonumber\hspace*{-30pt}\\[-8pt]\\[-8pt]
&&\qquad\le\Bigl|\sup_{y,\mathbf{x}}f_{Y|\X}(y|\mathbf{x})\Bigr|^2
\bigl[s(A)\bigr]2^p \int_{[-1,1]^p}\bigl|R_n(\t;
\tau,\mathbf{x})\bigr|^2\,d\t=O\bigl(h_n^{2s_2}
\bigr)\nonumber\hspace*{-30pt}
\end{eqnarray}
uniformly in $\tau\in[\delta^*,1-\delta^*]$ and $\mathbf{x}\in
\mathcal{D}$.

\item[\textit{Step} 4:] Now in view of assumption (A3), there exists $\lambda
_1>0$, such that
$g_n(\c,{\t}|\tau,\mathbf{x})\ge\lambda_1$ for all $\c$, $\t$
and $\mathbf{x}\in\mathcal{D}$
and $\tau\in[\delta^*,1-\delta^*]$, such that $|\c-\c_n(\mathbf
{x};\tau)|\le M_1$
and $(\c-\c_n(\mathbf{x};\tau))^\top\t(A)\ne0$. Let $\lambda_2$
be the smallest e-value of the
$s(A)\times s(A)$ matrix $\Gamma$. Then
for the first term on the right-hand side of (\ref{a11}), we have
%
%eA.23 #&#
\begin{eqnarray}\label{a12}
&& \biggl| \biggl[\int_{[-1,1]^p}g_n(\c,{\t}|\mathbf{x},\tau )
\t(A)\bigl\{\t(A)\bigr\}^\top \,d\mathbf{t} \biggr]\bigl(\c-\c_n(
\mathbf{x};\tau)\bigr) \biggr|
\nonumber\\[-8pt]\\[-8pt]
&&\qquad \ge \lambda_1\lambda_2\bigl|\c-
\c_n(\mathbf{x};\tau)\bigr|,\nonumber
\end{eqnarray}
\end{longlist}
for all $\c\in R^{s(A)}$ such that $|\c-\c_n(\mathbf{x};\tau)|\le M_1$.
The assertion in the proposition thus follows from
(\ref{a11}), (\ref{a13}), (\ref{a12}) and the conclusion reached in
step~1.
\end{pf}

We collect here some useful results for the verification of Euclidean
property of a class of functions.
\begin{longlist}[(C3)]
\item[(C1)]
Let $\mathfrak{F}=\{f(\cdot,t)\dvtx  t\in T\}$ be a class of functions indexed
by a
bounded subset $T $ of $R^d$. If there exists an $\alpha>0$ and a
nonnegative function
$\phi(\cdot)$ such that
\begin{eqnarray*}
&&\bigl|f(\cdot,t)-f\bigl(\cdot,t'\bigr)\bigr|\le\phi(\cdot)
\bigl\|t-t'\bigr\|^{\alpha}\qquad\mbox{for any }t, t'\in
T,
\end{eqnarray*}
then $\mathfrak{F}$ is Euclidean for the envelope $|f(\cdot,t_0)|+M \phi
(\cdot)$, where
$t_0$ is an arbitrary point of $T$ and $M=(2\sqrt{d}\sup_{T}\|
t-t_0\|)^{\alpha}$.
[Lemma 2.13 of \citet{r28}.]

\item[(C2)] If a class of functions $\mathfrak{F}$ is Euclidean for
an envelope $F$ and
$\mathfrak{g}$ is Euclidean for an envelope $G$, then $\{f+g\dvtx f\in
\mathfrak{F}, g\in\mathfrak{g}\}$
is Euclidean for the envelope $F+G$ and $\{fg\dvtx f\in\mathfrak{F}, g\in
\mathfrak{g}\}$
is Euclidean for the envelope $FG$. [Lemma 2.14 of \citet{r28}.]

\item[(C3)] Let $\lambda(\cdot)$ be a real-valued function of bounded
variation on $R$.
The class of all functions on $R^p$ of the form
$\{\lambda({\mathbf b}^\top\mathbf{x}+c)\dvtx  {\mathbf b}\in R^p, c\in R\}$
is Euclidean for a constant envelope. [Lemma 22(ii) of \citet{r27}.]

\item[(C4)] Let $\lambda(\cdot)$ be a real-valued function of bounded
variation on $R^+$.
The class of all functions on $R^p$ of the form
$\{\lambda(\|{\mathbf B} \mathbf{x}+{\mathbf b}\|)\dvtx  {\mathbf B}\in R^{m\times
p}, {\mathbf b}\in R^m\}$
is Euclidean for a constant envelope. [Lemma 22(i) of \citet{r27}.]
\end{longlist}

%
%coA.3 #&#
\begin{Corollary}\label{garden}
The following classes of functions are all Euclidean for an constant envelope:
$ \{I\{Y_{i}\le
Q_\tau(\X_{i})\}=I\{F(Y_{i}|\X_{i})\le
\tau\}, \tau\in(0,1)\}$,
$\{\X_{i\mathbf{x}}(h_{n},A)\dvtx  \mathbf{x}\in\mathcal{D}\}$,
$\{I(|\X_{i\mathbf{x}}|\le h_n)\dvtx  \mathbf{x}\in\mathcal{D}\}$ and
$\{I\{Y_{i}\le
Q_{n}(\X_{i},\mathbf{x};\tau)\}\dvtx  \mathbf{x}\in\mathcal{D},\tau
\in(0,1)\}$.
\end{Corollary}
\begin{pf} %{Proof of Corollary \ref{garden}}
This follows easily from (C2),
(C3) and (C4).
\end{pf}

\begin{pf*}{Proof of (\ref{lee})}
By Corollary \ref{garden},
any algebraic operations involving these classes of functions are also
Euclidean;
for example,  $\{\X_{ij}(h_{n},A)[I\{Y_{i}\le
Q_{n}(\X_{i},\X_j;\tau)\}-I\{Y_{i}\le
Q_\tau(\X_{i})\}]I(|\X_{ij}|\le h_n)\dvtx  \X_j\in\mathcal{D},\tau\in
(0,1)\}$.
This together with Theorem 37 in \citeauthor{r29} [(\citeyear{r29}), page~34] and the fact that
$Q_{n}(\X_{i},\X_j;\tau)-Q_\tau(\X_{i})=O(h_n^{s_2})$
lead to (\ref{lee}), that is,
with probability one,
\begin{eqnarray*}
&&\frac{1}{nh_n^p}\sum_{i}
\X_{ij}(h_{n},A)\bigl[I\bigl\{Y_{i}\le
Q_{n}(\X_{i},\X_j;\tau)\bigr\}-
I\bigl\{Y_{i}\le Q_\tau(\X_{i})\bigr\}
\bigr]I\bigl(|\X_{ij}|\le h_n\bigr)
\\
&&\qquad =o\bigl(n^{-1/2}\bigr)
\end{eqnarray*}
uniformly in $\X_j\in\mathcal{D}$, $\tau\in(0,1)$.
\end{pf*}
\end{appendix}

%
%
%
%accubitu suo, nardus mea dedit odorem suavitatis. Quoniam confortavit
%seras portarum tuarum, benedixit filiis tuis in te. Qui posuit fines
%tuos}

\section*{Acknowledgements} We thank the Editor, the Associate Editor
and two referees for
their helpful comments that have improved earlier versions of this paper.

%suskaldyti doi

% imsref loaded by linak, 2014-06-23 12:25:09
% imsref loaded by linak, 2014-06-23 12:51:04

\printaddresses
\end{document}